\documentclass{article}
\usepackage[english]{babel}
\usepackage[T1]{fontenc}
\usepackage[a4paper]{geometry}
\usepackage{amsfonts}
\usepackage{amsmath}
\usepackage{amsthm}
\usepackage{amssymb}
\usepackage{mathrsfs}
\usepackage{tikz-cd}
\usepackage{hyperref}
\usepackage{subcaption}
\usepackage{enumitem}
\usepackage{array}
\usepackage{mathtools,thmtools,thm-restate}
\usepackage{appendix}

\newtheorem{theorem}{Theorem}
\newtheorem*{theorem*}{Theorem}

\newtheorem{corollary}[theorem]{Corollary}
\newtheorem{proposition}[theorem]{Proposition}
\newtheorem{lemma}[theorem]{Lemma}
\theoremstyle{definition}
\newtheorem*{definition}{Definition}
\theoremstyle{remark}
\newtheorem*{remark}{Remark}

\renewcommand{\mod}{\ \mathrm{mod}\ }
\renewcommand{\O}{\mathcal{O}}
\newcommand{\C}{\mathbb{C}}
\newcommand{\R}{\mathbb{R}}
\newcommand{\N}{\mathbb{N}}

\DeclareMathOperator{\Id}{\text{Id}}

\newcommand{\HNC}{\mathbb{H}^{n}_{\C}}
\newcommand{\PHNC}{\partial\mathbb{H}^{n}_{\C}}

\newcommand{\PU}{\mathrm{PU}}
\newcommand{\U}{\mathrm{U}}
\DeclareMathOperator{\Mat}{\text{Mat}}

\newcommand{\PNC}{\C \mathbb{P}^{n}}

\DeclareMathOperator{\Alb}{\text{Alb}}
\newcommand{\HNL}{\mathbb{H}^{n}_{l}}
\newcommand{\PHNL}{\partial\mathbb{H}^{n}_{l}}
\DeclareMathOperator{\Ricci}{\text{Ricci}}
\DeclareMathOperator{\Stab}{\mathrm{Stab}}
\DeclareMathOperator{\Comm}{\mathrm{Comm}}

\title{Curvature properties and Shafarevich conjecture for toroidal compactifications of ball quotients}

\author{William Sarem}
\date{Octobre 13, 2025}

\begin{document} \maketitle

	\begin{abstract}
		We study toroidal compactifications of finite volume complex hyperbolic manifolds. We obtain results on the existence or nonexistence of Kähler metrics satisfying certain nonpositive curvature properties on these compactifications. Starting from quotients of complex hyperbolic space by deep enough non-uniform arithmetic lattices, we also verify the Shafarevich conjecture for their compactifications, by showing that their universal covers are Stein.
	\end{abstract}
	\section{Introduction}
	
	This article is a contribution to the study of closed Kähler manifolds obtained by \textit{toroidal compactifications} of spaces of the form $\HNC/\Gamma$, where $\HNC$ denotes the complex hyperbolic space of dimension $n$ and $\Gamma$ a ``deep enough'' torsion-free non-uniform lattice in $\PU(n,1)$. These compactifications, also called Mumford compactification, were first described in \cite{mumfordNewApproachCompactifying1975,ashSmoothCompactificationsLocally2010}. The metric study of these manifolds was done by Hummel and Schroeder \cite{hummelCuspClosingRank1996}. They showed:
	
	\begin{theorem*}[{\cite[Theorems 2 and 7]{hummelCuspClosingRank1996}}]\label{thm_HS}
		Let $\Gamma_0$ be a non-uniform torsion-free lattice in $\PU(n,1)$ for which the quotient space $\HNC/\Gamma_{0}$ admits a toroidal compactification $X_{\Gamma_0}$. Then $X_{\Gamma_0}$ is a Kähler manifold. Moreover, there exists a finite index subgroup $\Gamma'<\Gamma_0$ such that for any finite index subgroup $\Gamma<\Gamma'$ the toroidal compactification $X_{\Gamma}$ of $\HNC/\Gamma$ admits a Riemannian metric of nonpositive sectional curvature.
	\end{theorem*}
	
	The toroidal compactification of the quotient space $\HNC/\Gamma$ will be denoted by $X_\Gamma$, whenever it exists: precise conditions ensuring its existence are given in Section \ref{section_preli}, where we recall the construction of $X_{\Gamma}$. In fact, toroidal compactifications are projective manifolds \cite{mokProjectiveAlgebraicityMinimal2012,dicerboEffectiveResultsComplex2015}, and several authors have studied their properties, from the point of view of algebraic geometry  \cite{dicerboCanonicalDivisorSmooth2017,bakkerKodairaDimensionComplex2018,cadorelSymmetricDifferentialsComplex2021,memarianPositivityCotangentBundle2023}, or from other points of view  \cite{pyNoncoherentNonpositivelyCurved2017, dicerboClassificationArithmeticityToroidal2018-manuel}.
	
	In this article, we first show that these manifolds $X_\Gamma$ do not admit a Kähler metric with nonpositive sectional curvature (Theorem \ref{thm_1_pas_de_metrique}), confirming a claim from \cite{hummelCuspClosingRank1996}, but that they admit a Kähler metric with nonpositive holomorphic bisectional curvature (Theorem \ref{thm_2_metrique_HBC}). As a corollary, we answer in the affirmative a question raised by Diverio in \cite{diverioKobayashiHyperbolicityNegativity2022} on the existence of a closed complex manifold admitting a Kähler metric with quasi-negative holomorphic (bi)sectional curvature but admitting no Kähler metric with negative holomorphic sectional curvature (Corollary \ref{cor_repDiv}). In another direction, we also show that the universal cover $\widetilde{X_\Gamma}$ of $X_\Gamma$ is Stein when the lattice $\Gamma$ is arithmetic and deep enough.  To do so, we study the properties of the Albanese map of the manifolds obtained by compactifying finite covers of $\HNC/\Gamma$, and we obtain results concerning the Albanese image of the tori added during the compactification (Theorem \ref{thm_3_imm_alb}), extending earlier work of Eyssidieux \cite{eyssidieuxOrbifoldKahlerGroups2018}. This improvement allows us to prove that for deep enough finite index subgroups $\Gamma'$ of $\Gamma$, the universal cover of the toroidal compactification of $\HNC/\Gamma'$ is a Stein manifold (Theorem \ref{thm_4_rev_Stein}). When $n=2$ this result had been obtained earlier by Eyssidieux \cite{eyssidieuxOrbifoldKahlerGroups2018}.\newline

	Let us now describe more precisely and state these results.	In \cite{hummelCuspClosingRank1996}, Hummel and Schroeder assert that toroidal compactifications do not admit Kähler metrics with nonpositive sectional curvature. A proof of this assertion has been given in complex dimension 2 for ``deep enough'' lattices in \cite[Theorem A]{dicerboFinitevolumeComplexhyperbolicSurfaces2012}. Generalizing it to any dimension, we prove:
	
	\begin{theorem}\label{thm_1_pas_de_metrique}
		Let $\Gamma$ be a non-uniform torsion-free lattice in $\PU(n,1)$ for which the quotient space $\HNC/\Gamma$ admits a toroidal compactification $X_\Gamma$. Then there is no Kähler metric with nonpositive sectional curvature on $X_\Gamma$.
	\end{theorem}
	
	In particular, the techniques for constructing Kähler metrics on $X_\Gamma$ developed in \cite{hummelCuspClosingRank1996} cannot provide a metric with nonpositive sectional curvature. However, a computation shows that they are sufficient to construct a Kähler metric whose curvature verifies a weaker condition of nonpositivity, namely the nonpositivity of the holomorphic bisectional curvature:
	
	\begin{theorem}\label{thm_2_metrique_HBC}
		Let $\Gamma_0$ be a non-uniform torsion-free lattice in $\PU(n,1)$. There exists a finite index subgroup $\Gamma'<\Gamma_{0}$ such that for any finite index subgroup $\Gamma<\Gamma'$, the toroidal compactification $X_\Gamma$ of $\HNC/\Gamma$ admits a Kähler metric with nonpositive holomorphic bisectional curvature.
	\end{theorem}
	
	The metric constructed in the proof of this theorem has negative holomorphic bisectional curvature on the open subset $\HNC/\Gamma\subset X_{\Gamma}$. In particular, it is quasi-negative in the sense of \cite[Definition 1.1]{diverioQuasinegativeHolomorphicSectional2019}, i.e. the holomorphic bisectional curvature is nonpositive everywhere and negative at a point. In their article, Diverio and Trapani showed that every closed Kähler manifold with quasi-negative holomorphic sectional curvature has ample canonical bundle, a result which was already known under the assumption of negativity of the holomorphic sectional curvature \cite[Theorem 2]{wuNegativeHolomorphicCurvature2016} \cite[Corollary 1.3]{tosattiExtensionTheoremWu2017}. This answered a question of Yau. This situation naturally lead Diverio to ask in \cite{diverioKobayashiHyperbolicityNegativity2022} whether there exist closed Kähler manifolds of quasi-negative holomorphic sectional curvature but which do not admit a Kähler metric of negative holomorphic sectional curvature. Toroidal compactifications contain tori, which are not Kobayashi hyperbolic, so they allow to answer this question in the affirmative. Thus we obtain the following corollary.
	
	\begin{corollary}\label{cor_repDiv}
		For $\Gamma$ a deep enough non-uniform lattice as in Theorem \ref{thm_2_metrique_HBC}, $X_\Gamma$ is an example of a closed complex manifold which admits a Kähler metric of quasi-negative holomorphic bisectional curvature, but which does not admit a Kähler metric of negative holomorphic sectional curvature.
	\end{corollary}
	
	Along the way we observe that Theorem \ref{thm_2_metrique_HBC} combined with the main result of \cite{diverioQuasinegativeHolomorphicSectional2019} provides an alternative proof of \cite[Theorem 1.3]{dicerboCanonicalDivisorSmooth2017}: for $\Gamma$ a deep enough non-uniform lattice as in Theorem \ref{thm_2_metrique_HBC}, $X_\Gamma$ has ample canonical bundle. For a more direct proof, we remark in Section \ref{section_metrique_HBC} that the metric constructed in the proof of Theorem \ref{thm_2_metrique_HBC} has negative Ricci curvature. In addition, applying \cite[Theorem 3.1]{guenanciaQuasiprojectiveManifoldsNegative2022}, we get that any irreducible subvariety of $X_{\Gamma}$ which is not included in the complement of $\HNC/\Gamma$ is of general type.\newline
	
	We now state our results related to the Shafarevich conjecture for the compactifications $X_{\Gamma}$. This conjecture, stated in 1972 for projective manifolds, predicts that the universal cover of any closed Kähler manifold is holomorphically convex. It is known in the case where the fundamental group of the manifold is linear \cite{eyssidieuxLinearShafarevichConjecture2012,campanaRepresentationsLineairesGroupes2015}. 
	It is also known for specific examples \cite{eyssidieuxApplicationPeriodesVariation2015}. See also \cite{katzarkovNilpotentGroupsUniversal1997,katzarkovUniversalCoveringsAlgebraic1998,eyssidieuxConvexiteHolomorpheRevetements2004} for earlier results, and \cite{eyssidieuxLecturesShafarevichConjecture2011} for a survey article on this conjecture.\newline	
	
	For toroidal compactifications, if the lattice $\Gamma$ of $\PU(n,1)$ is deep enough, then by the theorem of Hummel and Schroeder recalled above, the universal cover $\widetilde{X_\Gamma}$ of $X_{\Gamma}$ is diffeomorphic to $\R^{2n}$ and can be endowed with a Kähler metric, thus it does not contain any compact analytic subvariety of positive dimension. Therefore, $\widetilde{X_\Gamma}$ is holomorphically convex if and only if it is Stein. When the lattice is arithmetic, Eyssidieux proved that, up to taking a finite index subgroup of the lattice, the Albanese map of $X_\Gamma$, denoted by $A_\Gamma\colon X_\Gamma \longrightarrow \Alb(X_\Gamma)$, is an immersion on the complement of the tori added during the compactification \cite[Theorem 1.3]{eyssidieuxOrbifoldKahlerGroups2018}. He then concluded in dimension 2 that the universal cover of $X_\Gamma$ is holomorphically convex, using a result of \cite{napierConvexityPropertiesCoverings1990} valid only in this dimension.	To study the higher dimensional case, we make Eyssidieux's result about Albanese mappings more precise. To do so, we analyse the behaviour of the Albanese map on the tori added during the compactification, which we will call in the following ``boundary tori of $X_\Gamma$''. More precisely, a boundary torus of $X_\Gamma$ is a connected component of the complement of $\HNC/\Gamma$ in $X_{\Gamma}$. With this terminology, we prove the following result.
	
	\begin{restatable}{theorem}{thmTrois}\label{thm_3_imm_alb}
		Let $\Gamma_0$ be a non-uniform torsion-free arithmetic lattice in $\PU(n,1)$. There exists a finite index subgroup $\Gamma' < \Gamma_0$ such that for any finite index subgroup $\Gamma < \Gamma'$, the toroidal compactification $X_{\Gamma}$ of $\HNC/\Gamma$ satisfies the following properties:\begin{itemize}
			\item Its Albanese map $A_\Gamma$ is an immersion on the open subset $\HNC/\Gamma \subset X_\Gamma$.
			\item For any boundary torus $T$ of $X_{\Gamma}$, the map $A_{\Gamma\,{\vert}T}$ is either an immersion or a constant.
		\end{itemize}
	\end{restatable}
	The arithmeticity of $\Gamma_{0}$ is used crucially in the proof of Theorem \ref{thm_3_imm_alb}, as in \cite{eyssidieuxOrbifoldKahlerGroups2018}. Deciding whether for every boundary torus $T$ of $X_{\Gamma_0}$, there exists a finite index subgroup $\Gamma<\Gamma_{0}$ for which the Albanese map $A_{\Gamma}$ of $X_\Gamma$ is an immersion when restricted to the tori above $T$ seems much more difficult. However the weaker result of Theorem \ref{thm_3_imm_alb} is enough for our purpose. We combine it with the arguments of \cite{napierConvexityPropertiesCoverings1990}, valid in any dimension, about holomorphic convexity with respect to a sufficiently positive line bundle, and we adapt the arguments given there in dimension 2 in order to generalize Eyssidieux's result to all dimension.
	
	\begin{theorem}\label{thm_4_rev_Stein}
		Let $\Gamma_0$ be a non-uniform torsion-free arithmetic lattice in $\PU(n,1)$ There exists a finite index subgroup $\Gamma' < \Gamma_0$ such that for any finite index subgroup $\Gamma < \Gamma'$, the universal cover $\widetilde{X_\Gamma}$ of the toroidal compactification $X_{\Gamma}$ of $\HNC/\Gamma$ is a Stein manifold.
	\end{theorem}
	
	This theorem is deduced from Proposition \ref{lem-PSH3} stated in Section \ref{section_rev_stein}, which is a more general result on projective manifolds which are the domain of a generically finite holomorphic map with values in a compact manifold with Stein universal cover.\newline
	
	The article is organized as follows. In Section \ref{section_preli}, we recall how toroidal compactifications of spaces of the form $\HNC/\Gamma$ are constructed. Then, Sections \ref{section_pas_de_metrique}, \ref{section_metrique_HBC}, \ref{section_imm_alb}, \ref{section_rev_stein} contain respectively the proofs of Theorems \ref{thm_1_pas_de_metrique}, \ref{thm_2_metrique_HBC}, \ref{thm_3_imm_alb}, \ref{thm_4_rev_Stein}. They use Section \ref{section_preli} and are independent of each other, except for Section \ref{section_rev_stein} which relies on Section \ref{section_imm_alb}.\newline

	\textbf{Acknowledgments.} I would like to express my deep gratitude to Pierre Py for his constant support and help throughout this work. This project was carried out during research internships at ETH Zürich and University of Rennes. I thank these institutions, and in particular Marc Burger and Christophe Dupont, for welcoming me. I also thank Simone Diverio and Matthew Stover for their valuable suggestions, as well as the anonymous referees for their comments and suggestions on the text. 
	
	\tableofcontents
	
	\section{Preliminaries}\label{section_preli}
	The complex hyperbolic space $\HNC$ is the unique complete simply connected Kähler manifold of dimension $n$ with constant holomorphic sectional curvature equal to $-4$. Its Riemannian sectional curvature is pinched between $-4$ and $-1$.  To construct it, we first consider a Hermitian form $\langle\cdot,\cdot\rangle$ of signature $(n,1)$ on $\C^{n+1}$. As a complex manifold, $\HNC$ is the projectivisation of the set of negative vectors for this form; it is an open subset of $\PNC$. Its boundary $\PHNC$ can be defined as its topological boundary in $\PNC$. We do not recall the construction of the hyperbolic metric on this space, and refer for example to \cite[Section 3.1]{goldmanComplexHyperbolicGeometry1999}. Its group of holomorphic isometries is the simple Lie group $\PU(n,1)$, it also acts on $\PHNC$. We will use the classical trichotomy of isometries into elliptic, parabolic or hyperbolic type, for which we refer for instance to  \cite[section II.6]{bridsonMetricSpacesNonpositive1999}.\newline
		
	We first introduce some notation which are analogous to the ones used in \cite{hummelCuspClosingRank1996,mokProjectiveAlgebraicityMinimal2012}. A nilpotent subgroup $N_{\xi}$ of $\PU(n,1)$ can be associated with every point $\xi \in \PHNC$, by taking the unipotent radical of the stabiliser of $\xi$ in $\PU(n,1)$. This Lie group is isomorphic to the Heisenberg group of dimension $2n-1$. The following calculations will be used in the proof of Proposition \ref{prop_toroidal_compactification}, they aim at providing global coordinates on the quotient of a horoball centered at $\xi$ by the action of a non-trivial element $Z$ in the center of $N_{\xi}$. \newline
		
	Let us fix a geodesic $\gamma$ parametrized with unit speed such that  $\lim_{t\to+\infty}\gamma(t) = \xi$. Denote also $o:=\gamma(0)$ and $\xi' := \lim_{t\to -\infty}\gamma(t)$. Then $\xi$ and $\xi'$ lift to vectors $f_1,f_2 \in \C^{n+1}$ which are isotropic for the Hermitian form defined above, non-collinear, and can be normalized so that $\langle f_1,f_2\rangle = 1$ and $o=[f_1-f_2]$. Here and in the remaining of the paragraph, $[\cdot]$ denotes the projection map $\C^{n+1}\setminus\{0\}\to \PNC$. Let $(f_3,\dots,f_{n+1})$ be a basis of $P :=(\C f_1\oplus \C f_2)^{\perp}$, orthonormal for the Hermitian product $\langle\cdot,\cdot \rangle \vert_{ P\times P}$. We identify $(P, \langle\cdot,\cdot \rangle \vert_{ P\times P})$ with $\C^{n-1}$ endowed with the standard Hermitian product whose associated norm is denoted by $\lVert\cdot\rVert$. In the basis $(f_1,\dots,f_{n+1})$ of $\C^{n+1}$, the form $\langle \cdot,\cdot\rangle$ has the following expression: \begin{equation*}
			\langle af_1+bf_2+v,af_1+bf_2+v\rangle = 2\Re(a\bar{b}) + \lVert v\rVert^{2}.
		\end{equation*}
		The group $N_{\xi}$ is made of the endomorphisms of $\C^{n+1}$ whose matrices in the basis $(f_1,\dots,f_{n+1})$ have the form\begin{align*}
			\left(\begin{array}{ccc}
				1 & -\frac{\lVert v\rVert^{2}}{2}-is & -{}^{t}\overline{v} \\
				0 & 1 & 0 \\
				0 & v & \mathrm{I}_{n-2}
			\end{array}\right),
		\end{align*}
		with $(s,v) \in\R \times \C^{n-1}$. We call $n_{(s,v)}$ the element of $\PU(n,1)$ defined by the matrix above. A global chart of $\HNC$, called the \textit{Siegel model }of $\HNC$, is given by
		\begin{align*}
			\widetilde{\Omega} := \{(a,v) \in \C\times\C^{n-1} \ \vert \ 2\Re(a)+\lVert v\rVert^{2}<0\} &\longrightarrow \HNC \\
			(a,v) &\longmapsto [af_1+f_2+v].
		\end{align*}
		Denote by $\widetilde{\mathcal{L}}$ the coordinates defined by the inverse of this map. For any real number $t$, let $a_t$ be the hyperbolic isometry with axis $\gamma$ defined in the basis $(f_1,\dots,f_{n+1})$ by the matrix $\mathrm{Diag}(e^{-t},e^{t},1,\dots,1)$. A calculation shows that the point $n_{(s,v)}\cdot a_t\cdot o$ has coordinates\begin{equation}\label{eq_coor_Ltilde}
			\widetilde{\mathcal{L}}\left(n_{(s,v)}\cdot a_t\cdot o\right) = \left(-\frac{\lVert v\rVert^{2}}{2}-is-e^{-2t},v\right).
		\end{equation}
		This shows in particular that the map $N_{\xi}\times \R \to \HNC$ which associates with each element $(n_{(s,v)},t)$ the point $n_{(s,v)}\cdot a_t\cdot o$ of the complex hyperbolic space is a diffeomorphism.\newline
		
		If $Z$ is an element of $N_{\xi}$ of the form $Z:=n_{(l,0)}$ with $l\ne 0$, the group $\langle Z\rangle$ which it generates acts freely and properly discontinuously on $\HNC$. The action of $Z$ is expressed in the coordinates $\widetilde{\mathcal{L}}$ by the translation $(a,u) \mapsto (a-il,u)$, thus the diagram of holomorphic maps
		\begin{equation}\label{diagramme}
			\begin{tikzcd}
				\HNC\arrow{d}{}\arrow{r}{\widetilde{\mathcal{L}}}&\widetilde{\Omega}\arrow{d}{(\exp(\frac{2\pi\bullet}{l}), \Id)} \\
				\HNC/\langle Z\rangle \arrow{r}{\mathcal{L}}  & \Omega.
			\end{tikzcd}
		\end{equation}
		commutes, where the map $\mathcal{L}$ defined by the diagram is a biholomorphism on its image\begin{equation}\label{eq-def-omega}
			\Omega := \{(b,v) \in \C\times \C^{n-1} \ \vert\ 0<\lvert b\rvert < \exp({\frac{-\pi \lVert v\rVert^{2}}{l}})\}.
		\end{equation}
		A horoball centered at $\xi$, denoted by $H_{\xi}$ or simply $H$ in what follows, is a subset of $\HNC \simeq N_{\xi}\times \R$ of the form $N_{\xi}\times (-\infty,t_0)$ for some real number $t_0$. This definition does not depend on the choice of the geodesic $\gamma$, for example because horoballs can also be defined as the sublevel sets of the Busemann function \cite[Section 4.1.2]{goldmanComplexHyperbolicGeometry1999}. The action of $Z$ on the hyperbolic space is identified with its action by left multiplication on the first factor of $N_{\xi}\times \R$. In particular, it preserves the horoballs. In coordinates $\widetilde{\mathcal{L}}$ and $\mathcal{L}$ respectively, a horoball $N_{\xi}\times (-\infty,t_0)$ and its quotient by the action of $Z$ identify with the open sets
		\begin{align*}
			\widetilde{\Omega}_{t_0} &:= \{(a,v) \in \C\times\C^{n-1} \ \vert\ \Re(a) < -\frac{\lVert v\rVert^{2}}{2}-\exp({-2t_0})\} \quad\text{and} \\
			\Omega_{t_0} &:= \{(b,v) \in \C\times\C^{n-1} \ \vert \ 0<\lvert b\rvert < \lambda(t_0)\exp({\frac{-\pi \lVert v\rVert^{2}}{l}})\},
		\end{align*}
		with $\lambda(t_0):= \exp(\frac{-2\pi e^{-2t_0}}{l})$.\newline
		
		A non-uniform lattice $\Gamma_0$ of $\PU(n,1)$ is a discrete subgroup of $\PU(n,1)$ for which the quotient space $\HNC/\Gamma_0$ has finite Riemannian volume, without being compact. Let $\Gamma_0$ be a non-uniform torsion-free lattice in $\PU(n,1)$, and $\Gamma$ a finite index subgroup of $\Gamma_0$. Then the quotient space $\HNC/\Gamma$ is a complex manifold naturally endowed with a Kähler metric of constant holomorphic sectional curvature $-4$. By Margulis' Lemma \cite{ballmannManifoldsNonpositiveCurvature1985} it admits a \textit{thick-thin decomposition }with a finite number of cusps, which means that we can write\begin{equation*}
			\HNC/\Gamma = Q \cup \coprod_{i=1}^{k}C_i,
		\end{equation*}
		where $Q$ is a compact subset of $\HNC/\Gamma$, $k\in \N^{*}$, and for any $1\le i\le k$, $C_i$ is a cusp, i.e., there exists a point $\xi_i$ of the boundary of $\HNC$ and a sufficiently deep horoball $H_i$ centered at $\xi_i$, globally fixed by the stabiliser $\Stab_{\Gamma}(\xi_i)$ of $\xi_i$ in $\Gamma$, such that $C_i$ is biholomorphically isometric to $H_i/\Stab_{\Gamma}(\xi_i)$ \cite{siuCompactificationNegativelyCurved1982}. The group $\Stab_{\Gamma}(\xi_i)$ contains no hyperbolic elements, thus it is a subgroup of $N_{\xi_i} \rtimes K_i$, where $K_i$ is the pointwise stabiliser in $\PU(n,1)$ of a geodesic ray ending at $\xi_i$. In particular, $\Stab_{\Gamma}(\xi_i)$ fixes each horosphere $N_{\xi_i} \times \{t\}$. Using the Auslander--Bieberbach theorem, one can show that there exists a finite subset $S\subset \Gamma_0$ such that for every finite index normal subgroup $\Gamma$ of $\Gamma_0$ which does not intersect $S$, the parabolic elements of $\Gamma$ are purely unipotent \cite{hummelRankOneLattices1998}. In this case, $\Stab_{\Gamma}(\xi_i)$ is a subgroup of $N_{\xi_i}$ for all $i$, and finally the cusp $C_i$ is isomorphic to the product of the nilmanifold $N_{\xi_i}/(\Gamma \cap N_{\xi_i})$ with the interval $(-\infty,t_0)$ for some real number $t_0$ \cite{hummelRankOneLattices1998}. Furthermore, since $H_i/(\Gamma \cap N_i)$ has finite volume, Fubini's theorem implies that $N_{\xi_i}/(\Gamma \cap N_{\xi_i})$ also has finite volume, and thus that $\Gamma \cap N_{\xi_i}$ is a lattice in $N_{\xi_i}$.
		
		We will say that $\Gamma$ is a lattice for which the quotient space $\HNC/\Gamma$ admits a toroidal compactification if $\Gamma$ is non-uniform, torsion-free and all parabolic subgroups of $\Gamma$ are purely unipotent, as explained in the previous paragraph. This is justified by the following proposition, which, for such lattices $\Gamma$, allows one to compactify the manifold $\HNC/\Gamma$ by identifying each cusp with a holomorphic punctured disk bundle over a complex torus of dimension $n-1$. The compactification can then be done by adding the zero section of this bundle. This identification is classical \cite[Section 2.1]{mokProjectiveAlgebraicityMinimal2012}, and we include it for the reader's convenience. 
		
		\begin{proposition} \label{prop_toroidal_compactification}
			Let $N$ be the $(2n-1)$-dimensional Heisenberg group, identified with $N_\xi$ for some point $\xi\in \PHNC$, $\Lambda$ a lattice in $N$, which acts on $\HNC \simeq N\times \R$ by left multiplication on the $N$-factor, and let $H:=N\times (-\infty,t_0)$ be a horoball of the complex hyperbolic space. Let $\rho:N\to \C^{n-1}$ be the surjective group morphism $n_{(s,v)}\mapsto v$. Then the map $\pi\colon H/\Lambda \to \C^{n-1}/\rho(\Lambda)$ is a holomorphic punctured disk bundle over the torus $T:=\C^{n-1}/\rho(\Lambda)$, which has negative curvature. More precisely, there exists a holomorphic line bundle $\widetilde{\pi}\colon L\to T$ on the torus $T$, endowed with a Hermitian metric $h$ with negative curvature $\Theta := -i\partial\bar{\partial}\log h$, as well as a holomorphic embedding $i:H/\Lambda \to L$ whose image is the open set\begin{equation*}
				\{v \in L \ \vert\ 0<\lVert v\rVert_h < 1\},
			\end{equation*}
			and such that $\pi = \widetilde{\pi}\circ i$.
		\end{proposition}
		
		\begin{proof}
			Since $\Lambda$ is a lattice in $N$, $\rho(\Lambda)$ is a lattice in $\C^{n-1}$ and in particular $T$ is a torus. The discrete subgroup $\Lambda \cap [N,N]$ is non-trivial, generated by an element $Z:=n_{(l,0)}$ for some $l>0$. 
			The map $\mathcal{L}$ defined above identifies $H/\langle Z\rangle$ with an open subset $\Omega_{t_0}$ of $\C\times \C^{n-1}$. Through this identification, $\Lambda$ acts on $\Omega_{t_0}$ by\begin{equation*}
				n_{(s,v)}\cdot (b,w) = \left(\exp\left({\frac{2\pi}{l}\left(-\frac{\lVert v\rVert^{2}}{2} - is - \langle w,v\rangle\right)}\right)b,v+w\right).
			\end{equation*}		
			The quotient space $H/\Lambda$ can thus be identified with the quotient of $\Omega_{t_0}$ by this action. Notice that this action naturally extends to an action of $\Lambda$ on $\C\times \C^{n-1}$. Let us endow the trivial line bundle map $\mathrm{pr}_2:\C\times \C^{n-1} \to \C^{n-1}$ with the Hermitian metric\begin{equation*}
				h((b,v),(b',v)) := \lambda(t_0)^{-2}\exp\left(\frac{2\pi \lVert v\rVert^{2}}{l}\right) b\bar{b'},
			\end{equation*}
			so that the inclusion map $j:\Omega_{t_0}\to \C\times\C^{n-1}$ defines a holomorphic embedding whose image is the open set\begin{equation*}
				\{(b,v) \in \C\times \C^{n-1} \ \vert\ 0<\lVert (b,v)\rVert_h < 1\}.
			\end{equation*}
			Since $\Omega_{t_0}$ is $\Lambda$-invariant, both the inclusion $j$ and the bundle map $\mathrm{pr}_2$ pass to the quotient, i.e. we have the following commutative diagram
			\begin{equation*}
				\begin{tikzcd}
					\Omega_{t_0}\arrow[hook]{d}{j}\arrow{r}{}&{\Omega}_{t_0}/\Lambda\arrow[hook]{d}{i} \\
					\C\times \C^{n-1} \arrow{d}{\mathrm{pr}_2}\arrow{r}{} & L \arrow{d}{\widetilde{\pi}} \\
					\C^{n-1}\arrow{r}&T.
				\end{tikzcd}
			\end{equation*}
			In this diagram, $L$ is the quotient of $\C\times\C^{n-1}$ by $\Lambda$, horizontal arrows are natural quotient maps and $\widetilde{\pi}$ is the only map making the bottom square of the diagram commute. Using the above expression for the action of $\Lambda$ on $\C\times \C^{n-1}$, we see that $\widetilde{\pi}:L\to T$ is a line bundle and $L$ can be endowed with the quotient Hermitian metric of $h$, still denoted by $h$, so that the image of $i$ is the open set
			\begin{equation*}
				\{v \in L \ \vert\ 0<\lVert v\rVert_h < 1\}.
			\end{equation*}
			The curvature of $(L,h)$, locally given by the form $-\frac{2\pi i}{l}\partial\bar\partial \lVert \cdot\rVert^{2}$, is negative.
		\end{proof}

		Let us denote by $X_{\Gamma}$ the toroidal compactification of $\HNC/\Gamma$ obtained by the previous construction. Recall that a line bundle $L\to X$ on a complex manifold $X$ is negative if it admits a Hermitian metric $h$ whose curvature form is negative, i.e. if the Chern class $c_1(L)$ of the bundle $L$ can be represented by a negative $(1,1)$-form. For each boundary torus $T$ added during the compactification, the normal bundle of $T$ in $X_{\Gamma}$ is identified with the normal bundle $(TL/TO)\vert_{ O}$ of the zero section $O$ in the total space of the bundle $L$, which is isomorphic to $L$. We deduce that: 
		
		\begin{corollary}\label{cor_fibre_normal_neg}
			For each boundary torus $T \subset X_{\Gamma}$ added during the compactification, the normal bundle of $T$ in $X_{\Gamma}$ is negative.
		\end{corollary}
		
		\section{Nonpositively curved Kähler manifolds containing submanifolds with $c_1=0$}\label{section_pas_de_metrique}
		
		In this section, we prove that for any non-uniform torsion-free lattice $\Gamma$ of $\PU(n,1)$, the toroidal compactification $X_\Gamma$, when defined, does not admit any Kähler metric with nonpositive sectional curvature (Theorem \ref{thm_1_pas_de_metrique}). We will use that $X_{\Gamma}$ contains a complex torus whose normal bundle is negative (Corollary \ref{cor_fibre_normal_neg}).
		The following proposition immediately implies the theorem. It uses the fact that the holomorphic sectional curvature decreases when passing to submanifolds, as well as relations between different notions of curvature, some of which are valid only for Kähler manifolds, and which are explained in \cite{diverioKobayashiHyperbolicityNegativity2022}.
		
		\begin{proposition}\label{prop_kahlerneg}
			Let $M$ be a complex manifold with a Kähler metric of nonpositive sectional curvature. Suppose that $M$ contains a compact submanifold $V$ whose first Chern class $c_1(TV)$ is zero. Then $V$ is totally geodesic and the Ricci curvature of $M$ is zero in restriction to $V$. Moreover the first Chern class $c_1(N_V)$ of the normal bundle of $V$ vanishes.
		\end{proposition}
		
		Before proving the proposition, we state our curvature conventions. On a Riemannian manifold $(M,g)$, we define:\begin{align}\label{conventionCourbure}
			&R(X,Y) := \nabla_{[X,Y]} - [\nabla_{X},\nabla_Y] \qquad \text{ and }\\
			&R(X,Y,Z,T) := g(R(X,Y)Z,T). \nonumber
		\end{align}
		Thus $(M,g)$ has nonpositive sectional curvature if and only  $R(X,Y,X,Y)\le 0$ for all vectors fields $X,Y$ of $M$. Accordingly, if $(f_1,\dots,f_n)$ is an orthonormal basis of $T_xM$, the Ricci curvature of $M$ at $x$ is defined  by\begin{equation*}
			\Ricci_M(X,Y) := \sum_{j=1}^{n}R(X,f_j,Y,f_j).
		\end{equation*}		
		With these conventions in mind, let us now turn to the proof of the Proposition.
		\begin{proof}[Proof of Proposition \ref{prop_kahlerneg}]
			The holomorphic sectional curvature of $V$ is less than or equal to that of $M$ restricted to the tangent bundle of $V$ \cite[chapter IX, proposition 9.2]{kobayashiFoundationsDifferentialGeometry1996}, therefore, it is nonpositive. By \cite[remark 2.8]{diverioKobayashiHyperbolicityNegativity2022}, in the case of a compact Kähler manifold, the non-positivity of the holomorphic sectional curvature implies the non-positivity of the scalar curvature, and the total scalar curvature (i.e. the integral of the scalar curvature) is zero if the Chern class $c_1(TV)$ of $V$ is zero. Therefore, the scalar curvature of $V$ is identically zero. Using \cite[proposition 2.9]{diverioKobayashiHyperbolicityNegativity2022}, we see that if at some point of $V$ there was a tangent vector with negative holomorphic sectional curvature, then the scalar curvature of $V$ would also be negative at this point, which is not the case. We deduce that the holomorphic sectional curvature of $V$ is identically zero. Because the latter determines the sectional curvature of $V$ \cite[chapter IX, proposition 7.1]{kobayashiFoundationsDifferentialGeometry1996}, we deduce that the metric $g$ restricted to $V$ is flat. The holomorphic sectional curvature of $V$ also bounds from below the holomorphic sectional curvature of $M$ restricted to the tangent bundle of $V$, so the latter is identically zero. From \cite[chapter IX, proposition 9.2]{kobayashiFoundationsDifferentialGeometry1996}, we deduce that the second fundamental form of $V$ is zero, thus that $V$ is totally geodesic.\newline
			
			Let us show that the Ricci curvature of $M$ restricted to $TV$ is identically zero. Let $k$ be the dimension of $V$, and $(z^{1},\dots,z^{n})$ be complex coordinates in the neighborhood of a point $x$ of $V$, such that $(\frac{\partial}{\partial z^{1}},\dots,\frac{\partial}{\partial z^{n}})$ forms an orthonormal basis of $T^{(1,0)}_xM$ and $(\frac{\partial}{\partial z^{1}},\dots,\frac{\partial}{\partial z^{k}})$ an orthonormal basis of $T^{(1,0)}_xV$. Then, writing $z^{i} =: x^{i} +\sqrt{-1}y^{i}$ and $e_i := \frac{\partial}{\partial x^{i}}$, we have that \begin{itemize}
				\item $(e_1, Je_1, \dots , e_n, Je_n)$ is an orthonormal basis of the tangent space $T_xM$ ;
				\item $(e_1, Je_1, \dots, e_{k}, Je_{k})$ is an orthonormal basis of the tangent space $T_xV$.
			\end{itemize}
			
			Since $e_1$ can be any unitary vector at $x$ tangent to $V$, it is enough to show that $\Ricci_M (e_1, e_1)$ vanishes. Recall that the first Bianchi identity gives for all vectors $X,Y \in T_xM$\begin{equation*}
				R(X,JX,Y,JY)	= R(X,Y,X,Y) + R(X,JY,X,JY).
			\end{equation*}
			In particular, if the sectional curvature is nonpositive, so is the holomorphic bisectional curvature. This formula also allows us to express the Ricci curvature at $x$:
			\begin{align*}
				\Ricci_M (e_1, e_1) &= \sum_{j=1}^n R(e_{1},e_{j},e_{1},e_{j}) + R(e_{1},Je_{j},e_{1},Je_{j})\\
				&= \sum_{j=1}^n R(e_1, Je_1, e_j, Je_j)\\
				&= \Ricci_V (e_1, e_1) + \sum_{j=k+1}^{n}R( e_1, Je_1,e_j, Je_j) \\
				&= \sum_{j=k+1}^{n}R( e_1, Je_1,e_j, Je_j).
			\end{align*}
			We assert that for all vectors $v,w$ tangent to $M$ at the same point, we have \begin{equation}\label{eq_courbure_determinant_polynome}
			R(v, Jv, w, Jw)^{2} \leq R(v, Jv, v, Jv) R(w, Jw, w, Jw).
			\end{equation}
			This inequality holds for all Kähler manifolds with nonpositive Riemannian sectional curvature. Provided that it is true, we apply it to $v=e_1$ and $w=e_j$ for all $j>k$: this shows that $\Ricci_{M\,\vert TV}=0$ since $R(e_1, Je_1, e_1, Je_1)=0$, due to the fact that the holomorphic sectional curvature of $M$ restricted to $TV$ is zero.\newline			
			
			The proof of Inequality \eqref{eq_courbure_determinant_polynome}, coming from \cite[Lemma 1]{mostowCompactKahlerSurface1980},\cite[Inequality (1)]{dicerboFinitevolumeComplexhyperbolicSurfaces2012}, consists in considering a one-parameter family of tangent planes at the considered point, and in expressing the curvature of these planes as a nonpositive function of the parameter. We use the fact that the metric is Kählerian to simplify the expression of the curvature. More precisely, for any real number $a$, we set
			\begin{align*}
				\left\{\begin{array}{l}
					p_a := av+w\\
					q_a := J(av-w).
				\end{array}
				\right.
			\end{align*}
			
			Since $p_{-a}=Jq_a$ and $q_{-a}=-Jp_a$, the map $a\mapsto R(p_a,q_a,p_a,q_a)$ is even. It can thus be written as a polynomial $P(a^{2})$ of degree at most 2 in $a^2$:
			\begin{equation*}
				P(a^{2}) = R(p_a,q_a,p_a,q_a) = R(v,Jv,v,Jv)a^4 + Aa^2 + R(w,Jw,w,Jw),
			\end{equation*}
			with\begin{align*}
				A := &R(v,Jv,w,-Jw) + R(v,-Jw,v,-Jw) + R(v,-Jw,w,Jv) +\\ &\qquad R(w,Jv,v,-Jw) + R(w,Jv,w,Jv) + R(w,-Jw,v,Jv) \\
				= &-R(v,Jv,w,Jw) + R(v,-Jw,v,-Jw) - R(v,-Jw,v,-Jw) +\\ & \qquad R(w,Jv,v,-Jw) - R(w,Jv,v,-Jw)-R(v,Jv,w,Jw) \\
				=&-2R(v,Jv,w,Jw).
			\end{align*}
			Thus\begin{equation*}
				P(a^{2}) = R(v,Jv,v,Jv)a^4 -2R(v,Jv,w,Jw) a^2 +  R(w,Jw,w,Jw).
			\end{equation*}
			\
			
			Since the sectional curvature is nonpositive, the polynomial $P$ is nonpositive on $\R_+$. Recall also that $R(X,JX,Y,JY)\le 0$ for all vectors $X,Y \in T_xM$. We distinguish two cases.\begin{itemize}
				\item If the coefficient $R(v,Jv,v,Jv)$ vanishes, then when $a$ becomes increasingly large,	 the inequality $R(p_a,q_a,p_a,q_a)\le 0$ forces  $R(v,Jv,w,Jw)$ to vanish as well.
				\item Otherwise, set\begin{equation*}
					a=\left(\frac{R(v,Jv,w,Jw)}{R(v,Jv,v,Jv)}\right)^{1/2}.
				\end{equation*} Writing the nonpositivity of $P(a)$ gives the desired inequality.
			\end{itemize} In all cases, we get Inequality \eqref{eq_courbure_determinant_polynome}.\newline
			
			Let us now show that the normal bundle $N_V$ of $V$ has vanishing first Chern class. Multiplicativity of the total Chern classes with respect to the exact sequence\begin{equation*}
				0 \longrightarrow TV \longrightarrow TM\vert_{V} \longrightarrow N_V \longrightarrow 0
			\end{equation*}
			yields the formula $c_1(N_V) = c_1(TM)\vert_{V} - c_1(TV)$.
			Moreover, if $\rho := \Ricci_M(\cdot,J\cdot)$ denotes the Ricci form of $M$, the class $c_1(TM)$ is represented by $-\frac{\sqrt{-1}}{2\pi}\rho$ \cite[Chapter 2]{kobayashiDifferentialGeometryComplex1987}, which vanishes in restriction to $V$. Since $c_1(TV)=0$ by hypothesis, we deduce that $c_1(N_V)=0$.
		\end{proof}

		\section{A metric on $X_{\Gamma}$ with nonpositive holomorphic bisectional curvature}\label{section_metrique_HBC}

		In this section, we prove that for any non-uniform torsion-free lattice $\Gamma_0$ of $\PU(n,1)$, there exists a finite index subgroup $\Gamma'<\Gamma_{0}$ such that for any finite index subgroup $\Gamma<\Gamma'$, the toroidal compactification $X_\Gamma$ of $\HNC/\Gamma$ admits a Kähler metric with nonpositive holomorphic bisectional curvature (Theorem \ref{thm_2_metrique_HBC}). We will need the following lemma, whose proof we omit.\begin{lemma}\label{lem-cutoff}
			There exists a constant $C_0>0$ such that for every $C\ge C_0$, there is a smooth function $f:[0,C]\to \R_+$ which coincides with $\cosh$ on some neighborhood of $0$, coincides with $\exp$  on some neighborhood of $C$ and such that the functions $f,f',f'',f'''$ are positive on $(0,C]$.
		\end{lemma}
		
		The proof of Theorem \ref{thm_2_metrique_HBC} takes up Hummel and Schroeder's construction of a Kähler metric on toroidal compactifications \cite{hummelCuspClosingRank1996}, passing to finite index subgroups in order to get space for the construction of a metric with nonpositive holomorphic bisectional curvature. In a nutshell, the idea is to work separately in each cusp, and, provided that the cusp is ``big enough'', to patch there some compact manifold, endowed with a Kähler metric which glues to the complex hyperbolic one, and whose bisectional curvature is nonpositive. For the reader's convenience, we first provide a detailed description of the construction of the metric, (steps 1 to 5) as carried out by Hummel and Schroeder, and then we come to the original part of the proof, which is the computation of the holomorphic bisectional curvature (steps 6 and 7).
		\newline
		
		\underline{Step 1.} Notation and description of a fixed cusp.
		\vspace{.5\baselineskip}
		
		Using the model $(s,v) \mapsto n_{(s,v)}$ of the Heisenberg group $N$ described in Section \ref{section_preli}, we endow $N$ with the left-invariant metric $\mu$ expressed at the identity point by\begin{equation*}
				\mu_{\Id}:=ds^{2}\oplus dv^{2}.
			\end{equation*}Let $\mathfrak{n}$ be the Lie algebra of $N$, $Z$ the generator of $[\mathfrak{n},\mathfrak{n}]\simeq \R$ such that $\exp(Z)=n_{(1,0)}$ and $\mathfrak{r}$ the $\mu$-orthogonal of $[\mathfrak{n},\mathfrak{n}]$. Assuming that $\Gamma_0$ is a torsion-free non-uniform lattice in $\PU(n,1)$ with purely unipotent parabolic elements  as explained in Section \ref{section_preli}, any cusp in $\HNC/\Gamma_0$ is biholomorphically isometric to\begin{equation*}
				((-\infty,A_0]\times N/\Lambda_0, \mu_{0}, J_0),
			\end{equation*}
			where $A_0$ is some real number and $\Lambda_0$ a lattice in $N$ which depend on the cusp; $\mu_0$ is the $N$-invariant metric expressed at a point $(t,\Id\! \mod \Lambda_0)$ by \begin{equation*}
				\mu_{0 (t,\Id\! \mod \Lambda_0) }:=dt^{2}\oplus \exp({2t})\mu\vert_{ \mathfrak{r}}\oplus \exp({4t})\mu\vert_{ \R Z};
			\end{equation*}and $J_0$ is the $N$-invariant complex structure defined at a point $(t,\Id \mod \Lambda_0)$ by\begin{equation*}\left\{
			\begin{array}{l}
				J_0Z =\exp({2t})\frac{\partial}{\partial t},\\
				\forall v\in \C^{n-1}, \exp\circ J_0\circ\exp^{-1}(n_{(0,v)})=n_{(0,iv)}.
\end{array}
\right.
\end{equation*}
\newline

\underline{Step 2.} Passing to a finite index subgroup of the lattice.
\vspace{.5\baselineskip}

This is achieved in the following lemma, which is a direct consequence of Malcev's theorem. 

\begin{lemma}\label{lem-thm4-ssgrp}
	Let $\Gamma_0$ be a non-uniform torsion-free lattice in $\PU(n,1)$, and $C_0$ be the constant given by Lemma \ref{lem-cutoff}. Then there exists a finite index subgroup $\Gamma'<\Gamma_{0}$ such that for any finite index subgroup $\Gamma<\Gamma'$, any cusp of $\HNC/\Gamma$, which lies above some cusp $((-\infty,A_0]\times N/\Lambda_0, \mu_{0}, J_0)$ of $\HNC/\Gamma_{0}$, is biholomorphically isometric to $((-\infty,A_0]\times N/\Lambda, \mu_{0}, J_0)$ for some finite index subgroup $\Lambda<\Lambda_0$ such that the group $\Lambda\cap [N,N]$ is generated by $n_{(l,0)}$ with $l>2\pi \exp({2(C_0-A_0)})$.
\end{lemma}
\

\underline{Step 3.} Description of Hummel and Schroeder's model.
\vspace{.5\baselineskip}

We now fix $\Gamma_0$ a non-uniform torsion-free lattice in $\PU(n,1)$ and $\Gamma$ a finite index subgroup of $\Gamma_{0}$ as in Lemma \ref{lem-thm4-ssgrp}. Let $((-\infty,A_0]\times N/\Lambda, \mu_{0}, J_0)$ be a cusp of $\HNC/\Gamma$, and $l>2\pi \exp({2(C_0-A_0)})$ be the positive number such that $\Lambda\cap [N,N]$ is generated by $n_{(l,0)}$.
Let $\varphi$ be the automorphism of $N$ given by\begin{equation*}
	\varphi(n_{(s,v)}) = n_{\left({\frac{2\pi}{l}}s,\sqrt\frac{2\pi}{l}v\right)}.
\end{equation*}
The next lemma describes the model used by Hummel and Schroeder to show that $X_{\Gamma}$ is a Kähler manifold, and in which we will define a Kähler metric of nonpositive holomorphic bisectional curvature.

\begin{lemma}
			 Set $C:=\ln\sqrt{\frac{l}{2\pi}}$, $A:=A_0+C$, and let $f:[0,A]\to \R_+$ be the function given by Lemma \ref{lem-cutoff} and $g:=ff'$. Then $((-\infty,A_0]\times N/\Lambda, J_0)$ is biholomorphic to $((0,A]\times N/\varphi(\Lambda),J)$, where $J$ is the $N$-invariant complex structure defined at a point $(t,\Id\!\mod \varphi(\Lambda))$ by\begin{equation*}\left\{
				\begin{array}{l}
					JZ =g(t)\frac{\partial}{\partial t},\\
					\forall v\in \C^{n-1}, \exp\circ J\circ\exp^{-1}(n_{(0,v)})=n_{(0,iv)}.
				\end{array}
				\right.
			\end{equation*}
\end{lemma}
		
\begin{proof}
	First define the diffeomorphism\begin{equation*}
		\Phi :\left\{\begin{array}{rcl}
			(-\infty,A_0]\times N/\Lambda &\longrightarrow &(-\infty,A]\times N/\varphi(\Lambda) \\
			(t,n \mod \Lambda) &\longmapsto &(t+C,\varphi(n) \mod \varphi(\Lambda)).
		\end{array}
		\right. 
	\end{equation*}
	Easy computations show that $\Phi_*\mu_0 = \mu_0$ and $\Phi_*J_0=J_0$. In this new model, we have by construction  $\varphi(n_{(l,0)})=n_{(2\pi,0)}$ and $A> C_0$.
	
	Then let $\psi:(0,A] \to (-\infty,A]$ be the solution of the following differential equation\begin{equation*}
		\left\{\begin{array}{rl}
			\psi'(t)=&\frac{\exp({2\psi(t)})}{g(t)} \\
			\psi(A)=&A
		\end{array}
		\right.
	\end{equation*}
	It is easily seen that $\psi$ is a well defined diffeomorphism, and that $\psi$ is the identity in a neighborhood of $A$. Define \begin{equation*}
		\Psi:=\psi\times \Id\colon(0,A]\times N/\varphi(\Lambda) \longrightarrow (-\infty,A]\times N/\varphi(\Lambda).
	\end{equation*}
	A computation shows that $J:=\Psi^{*}J_0$ is the $N$-invariant complex structure described in the Lemma. Thus $\Psi$ provides the required biholomorphism.
\end{proof}
	
\underline{Step 4.} A Kähler metric on the cusp.
\vspace{.5\baselineskip}

We endow the complex manifold $((0,A]\times N/\varphi(\Lambda),J)$ with the $N$-invariant Hermitian metric $\mu_{f,g}$ expressed at a point $(t,\Id \mod \varphi(\Lambda))$ by \begin{equation*}
	\mu_{f,g \ (t,\Id\! \mod \varphi(\Lambda)) }:=dt^{2}\oplus f(t)^{2}\mu\vert_{ \mathfrak{r}}\oplus g(t)^{2}\mu\vert_{ \R Z }.
\end{equation*}
		
Notice that the metric $\mu_{f,g}$ coincides with $\Psi^{*}\mu_0$ on a neighborhood of the boundary of the manifold, hence it glues back to the hyperbolic metric on the thick part of $\HNC/\Gamma$. We are going to use the following facts from \cite{hummelCuspClosingRank1996}.
		
		\begin{lemma}
			\begin{enumerate}
			\item The toroidal compactification of the cusp can be realized in this model in the following way: using polar coordinates, we first identify $(0,A)\times N/\langle n_{(2\pi,0)}\rangle$ with $\mathbb{D}^{*}(0,A)\times \C^{n-1}$, where $\mathbb{D}^{*}(0,A)$ is the punctured disk in $\C$ with radius $A$; then we take the quotient by the action of $\varphi(\Lambda)$ \cite[Proof of Theorem 7]{hummelCuspClosingRank1996}.
			\item The equality $g=ff'$ implies that the metric $\mu_{f,g}$ is Kählerian \cite[Proof of Theorem 7]{hummelCuspClosingRank1996}.
			\item The facts that $f=\cosh$ and that $g=ff'$ imply that $\mu_{f,g}$ extends to the compactification of the cusp \cite[Lemma 3.8]{hummelCuspClosingRank1996}.
		\end{enumerate}
		\end{lemma}
		
\underline{Step 5.} Some preliminary curvature computations.
\vspace{.5\baselineskip}

			As in \cite{hummelCuspClosingRank1996}, we use the curvature convention  given by Equation  \eqref{conventionCourbure} and we see the Riemannian curvature $R$ of $\mu_{f,g}$ as a symmetric form defined on bivectors. In the next lemma, which will be used to compute the holomorphic sectional curvature of $\mu_{f,g}$, we identify the tangent space of $(0,A) \times N/\varphi(\Lambda) $ at a point of the form $(t,\Id \mod \varphi(\Lambda))$ with $\R\oplus \mathfrak{n}$.
			
			\begin{lemma}[{\cite[Appendix]{hummelCuspClosingRank1996}}]\label{lem-thm4-calculHS}
				Let $X\in \mathfrak{r}$ be a tangent vector at $(t,\Id \mod \varphi(\Lambda))$ which satisfies $\mu(X,X)=1$. Define the bivectors\begin{align*}
					(E_1,E_2,E_3,E_4,E_5,E_6) := \biggl(
					&\frac{X \wedge JX}{\lVert X\wedge JX \rVert},
					\frac{Z \wedge JZ}{\lVert Z\wedge JZ\rVert},
					\frac{X \wedge Z}{\lVert X\wedge Z\rVert}, \\
					&\left.\frac{JX \wedge JZ}{\lVert JX\wedge JZ\rVert},
					\frac{JX \wedge Z}{\lVert JX\wedge Z\rVert},
					\frac{X \wedge JZ}{\lVert X\wedge JZ\rVert}
					\right).
				\end{align*}
				Then \begin{equation*}
					(R(E_i,E_j))_{1\le i,j\le 6} =: \left(\begin{array}{ccc}
						G & 0 & 0 \\
						0 & F & 0 \\
						0 & 0 & F
					\end{array}\right),
				\end{equation*}
				where\begin{align*}
					G = \left(\begin{array}{cc}
						-4\left(\frac{f'}{f}\right)^{2} & -2\frac{f''}{f} \\
						-2\frac{f''}{f} & -3\frac{f''}{f} - \frac{f'''}{f'}
					\end{array}\right) =: \left(\begin{array}{cc}
						-h_1^{2} & -2k \\
						-2k & -h_2^{2}
					\end{array}\right), \\
					F = \left(\begin{array}{cc}
						-\frac{f''}{f} & -\frac{f''}{f} \\
						-\frac{f''}{f} & -\frac{f''}{f}
					\end{array}\right) =: \left(\begin{array}{cc}
						-h_3^{2} & -k \\
						-k & -h_4^{2}
					\end{array}\right).
				\end{align*}
				The functions $h_i:(0,A)\to (0,+\infty)$ defined through $G$ and $F$ are well defined because of the positivity conditions for $f$ and its derivatives. Here as in the sequel, the functions $f,g,k,h_i$ shall implicitly be evaluated at $t$.
			\end{lemma}
						
\underline{Step 6.} Computation of the holomorphic bisectional curvature of $\mu_{f,g}$.
			
This is done in the following proposition.
			\begin{proposition}\label{prop-thm4-HBC}
				Let $Y, \Xi$ be two tangent vectors at $(t,\Id \mod \varphi(\Lambda))$ for some $t\in (0,A)$. Write $Y=:aX+bZ+cJZ$ and $\Xi =: \alpha \widetilde{X} + \beta Z + \gamma JZ$ with $a,b,c,\alpha,\beta,\gamma \in \R$, and $X,\widetilde{X} \in \mathfrak{r}$ of norm $\mu(X,X)=\mu(\widetilde{X},\widetilde{X})=1$. Also write $\widetilde{X}=dX+eJX+W$ with $d,e\in \R$ and $W$ a vector $\mu$-orthogonal to $X$ and $JX$.
				For $\mu_{f,g}$, the holomorphic bisectional curvature $R(Y \wedge JY,\Xi \wedge J\Xi)$ of $Y$ and $\Xi$ is then given by:
				\begin{align}\label{form-BHC}
					-R(&Y \wedge JY,\Xi \wedge J\Xi) = a^{2}\alpha^{2}\left((d^{2}+e^{2})f^{4}h_1^{2} + \frac{2g^{2}}{f^{2}}\lVert W\rVert^{2}\right) + \nonumber\\ & h_2^{2}g^{4}(b^{2}+c^{2})(\beta^{2}+\gamma^{2}) + 2f^{2}g^{2}h_3^{2}
					\big(a^{2}(\beta^{2}+\gamma^{2}) +\alpha^{2}(b^{2}+c^{2})\ + \\
					&2a\alpha(b\beta + c\gamma)\mu(X,\widetilde{X})
					+ 2a\alpha(c\beta - b\gamma)\mu(X,J\widetilde{X})\big).\nonumber
				\end{align}
			\end{proposition}

			\begin{proof}
				We break down the proof of Proposition \ref{prop-thm4-HBC} into 5 steps.
				
				{\parindent30pt\textit{Step 1.} Decomposition of $Y \wedge JY$ and $\Xi \wedge J\Xi$.}
					\vspace{.5\baselineskip}
					
				We compute that\begin{align*}
				Y \wedge JY = & B_1 + B_2 + B_3 \qquad \text{ with: }\\
				&B_1 = a^{2}X\wedge JX + (b^{2}+c^{2})Z\wedge JZ, \\	
				&B_2 = -ac(X\wedge Z + JX\wedge JZ),\\
				&B_3 =-ab(JX\wedge Z + (-X\wedge JZ)).
			\end{align*}
			In the same way, $\Xi \wedge J\Xi = \beta_1 + \beta_2 + \beta_3$ with $\beta_i$ defined as $B_i$ by replacing $a,b,c$ by $\alpha,\beta,\gamma$ and $X$ by $\widetilde{X}$.
			\newline
			
				{\parindent30pt\textit{Step 2.} Computation of $\mathcal{R} := R(X\wedge JX,\widetilde{X}\wedge
				 J\widetilde{X})$.}
				\vspace{.5\baselineskip}
					 
				\noindent We compute that
				 \begin{align*}
				 	\mathcal{R} &= R(X,JX,dX+eJX+W,-eX+dJX+JW)\\
				 	&=(d^{2}+e^{2})R(X,JX,X,JX) + dR(X,JX,X,JW)\ +\\ &\qquad\qquad dR(X,JX,W,JX) +eR(X,JX,JX,JW)\ -\\
				 	&\qquad\qquad eR(X,JX,W,X) + R(X,JX,W,JW) \\
				 	&= (d^{2}+e^{2})R(X,JX,X,JX)+2R(X,JX,X,eW+dJW)\ +\\ &\qquad\qquad R(X,JX,W,JW) \\
				 	&= -(d^{2}+e^{2})f^{4}h_1^{2} -2(\frac{3g^{2}}{f^{2}}+f'^{2})\langle JX,eW+dJW\rangle\ + \\&\qquad\qquad R(X,JX,W,JW)\\
				 	&= -(d^{2}+e^{2})f^{4}h_1^{2} + 0 + R(X,JX,W,JW).
				 \end{align*}
				 To pass from the third to the fourth line, we used Lemma \ref{lem-thm4-calculHS} and the computation of $R(X,JX)X$ done in \cite[page 305]{hummelCuspClosingRank1996}. Moreover,\begin{align*}
				 	R(X,JX,W,JW)&=\langle \nabla_{[X,JX]}W-\nabla_X\nabla_{JX}W + \nabla_{JX}\nabla_X W,JW \rangle \\
				 	&= 2\langle \nabla_ZW,JW \rangle \qquad \text{ because } \nabla_X W = \nabla_{JX} W = 0\\ &\hspace{5.4cm} \text{ and }[X,JX]=2Z\\
				 	&=\frac{-2g^{2}}{f^{2}}\lVert W\rVert^{2} \qquad \text{ from \cite[(3.2)]{hummelCuspClosingRank1996}.}
				 \end{align*}
				 In conclusion\begin{equation*}
				 	-R(X\wedge JX,\widetilde{X}\wedge J\widetilde{X}) = (d^{2}+e^{2})f^{4}h_1^{2} + \frac{2g^{2}}{f^{2}}\lVert W\rVert^{2}.
				 \end{equation*}

				{\parindent30pt\textit{Step 3.} Computation of $R(B_1,\beta_1)$.}
					\vspace{.5\baselineskip}
					
				We compute that\begin{align*}
					R(B_1,\beta_1) &= a^{2}\alpha^{2}R(X\wedge JX,\widetilde{X}\wedge J\widetilde{X})\ -\\ &\qquad2k a^{2}(\beta^{2}+\gamma^{2})\lVert X\wedge JX\rVert \lVert Z\wedge JZ\rVert\ - \\ 
					&\qquad  2k \alpha^{2}(b^{2}+c^{2})\lVert \widetilde{X}\wedge J\widetilde{X}\rVert \lVert Z\wedge JZ\rVert\ - \\
					&\qquad h_2^{2}(b^{2}+c^{2})(\beta^{2}+\gamma^{2})\lVert Z\wedge JZ\rVert^{2}\\
					&= a^{2}\alpha^{2}R(X\wedge JX,\widetilde{X}\wedge J\widetilde{X})\ -\\
					&\qquad2h_3^{2}f^{2}g^{2}\left(a^{2}(\beta^{2}+\gamma^{2}) +\alpha^{2}(b^{2}+c^{2}) \right)\ -\\ &\qquad h_2^{2}g^{4}(b^{2}+c^{2})(\beta^{2}+\gamma^{2})\\
					&=-a^{2}\alpha^{2}\left((d^{2}+e^{2})f^{4}h_1^{2} + \frac{2g^{2}}{f^{2}}\lVert W\rVert^{2}\right)\ - \\ &\qquad 2h_3^{2}f^{2}g^{2}\left(a^{2}(\beta^{2}+\gamma^{2}) +\alpha^{2}(b^{2}+c^{2}) \right)\ -\\ &\qquad h_2^{2}g^{4}(b^{2}+c^{2})(\beta^{2}+\gamma^{2}).
				\end{align*}
				
				{\parindent30pt\textit{Step 4.} Computation of $R(B_i,\beta_j)$ for $i,j\in \{2,3\}$.}
					\vspace{.5\baselineskip}
					
				First notice that the map $(U,\widetilde{U}) \mapsto R(U\wedge Z,\widetilde{U}\wedge Z)$ is bilinear symmetric, of quadratic form\begin{equation*}
					R(U\wedge Z,U\wedge Z) = -f^{2}g^{2}h_3^{2}\mu(U,U).
				\end{equation*}
				By splitting, we deduce that for all $U,\widetilde{U}$, we have \begin{equation*}
					R(U\wedge Z,\widetilde{U}\wedge Z) = -f^{2}g^{2}h_3^{2}\mu(U,\widetilde{U}).
				\end{equation*}
				Then, using that the metric is Kählerian, we find that\begin{align*}
					R(B_2,\beta_2) &= ac\alpha\gamma R(X\wedge Z + JX \wedge JZ,\widetilde{X}\wedge Z + J\widetilde{X} \wedge JZ) \\
					&= 4ac\alpha\gamma R(X\wedge Z,\widetilde{X}\wedge Z) \\
					&= -4f^{2}g^{2}h_3^{2}\mu(X,\widetilde{X})ac\alpha\gamma.\\
					R(B_3,\beta_3) &= 4ab\alpha\beta R(JX\wedge Z,J\widetilde{X}\wedge Z)\\
					&= -4f^{2}g^{2}h_3^{2}\mu(JX,J\widetilde{X}) ab\alpha\beta\\
					&= -4f^{2}g^{2}h_3^{2}\mu(X,\widetilde{X}) ab\alpha\beta.\\
					R(B_2,\beta_3) &= 4ac\alpha\beta R(X\wedge Z,J\widetilde{X}\wedge Z)\\
					&=-4f^{2}g^{2}h_3^{2}\mu(X,J\widetilde{X})ac\alpha\beta.\\
					R(B_3,\beta_2) &= 4ab\alpha\gamma R(JX\wedge Z,\widetilde{X}\wedge Z)\\
					&=-4f^{2}g^{2}h_3^{2}\mu(JX,\widetilde{X})ab\alpha\gamma\\
					&=-4f^{2}g^{2}h_3^{2}\mu(X,J\widetilde{X})ab\alpha\gamma.
				\end{align*}
				
				{\parindent30pt\textit{Step 5.} Vanishing of the remaining terms.}
					\vspace{.5\baselineskip}
					
			According to the computations of the appendix of \cite{hummelCuspClosingRank1996}, $R(X,JX)Z$ is a multiple of \linebreak $JZ=g\frac{\partial}{\partial t}$, which is orthogonal to both $\widetilde{X}$ and $J\widetilde{X}$. Therefore the two terms $R(B_1,\beta_2)$ and $R(B_1,\beta_3)$ vanish. By symmetry, we also have $R(B_2,\beta_1) = R(B_3,\beta_1) =0$.
			
			Putting \textit{Steps} 2 to 5 together, we obtain Formula \eqref{form-BHC}.
			\end{proof}

\underline{Step 7.} Conclusion: nonpositivity of the holomorphic bisectional curvature.
	\vspace{.5\baselineskip}
	
To show that $R(Y \wedge JY,\Xi \wedge J\Xi)\le0$, it is enough to show that the term\begin{equation*}
	a^{2}(\beta^{2}+\gamma^{2}) +\alpha^{2}(b^{2}+c^{2}) +
	2a\alpha(b\beta + c\gamma)\mu(X,\widetilde{X})
	+ 2a\alpha(c\beta - b\gamma)\mu(X,J\widetilde{X})
\end{equation*}
in Formula \eqref{form-BHC} is nonpositive. This is a consequence of the following lemma.

			\begin{lemma}
				\begin{equation*}
					\left\lvert 2a\alpha\left((b\beta + c\gamma)\mu(X,\widetilde{X})
					+ (c\beta - b\gamma)\mu(X,J\widetilde{X})\right) \right\rvert \le a^{2}(\beta^{2}+\gamma^{2}) +\alpha^{2}(b^{2}+c^{2}).
				\end{equation*}
			\end{lemma}
			
			\begin{proof}
				The two real numbers $x:= \mu(X,\widetilde{X})$ and $y:= \mu(X,J\widetilde{X})$ verify the inequality\begin{equation*}
				x^{2}+y^{2} = \lvert x+iy\rvert^{2} = \lvert \mu_{\C}( X,\widetilde{X})\rvert^{2} \le  \mu(X,X)\mu(\widetilde{X},\widetilde{X}) \le 1,
			\end{equation*}
			where $\mu_{\C}$ is the hermitian product associated with $\mu\vert_{ \mathfrak{r}}$, defined by \begin{equation*}
				\mu_{\C}(U,\widetilde{U}):= \mu(U,\widetilde{U})+i\mu(U,J\widetilde{U}).
			\end{equation*}
			Using the triangle inequality and the inequality $2XY\le X^{2}+Y^{2}$, we get that\begin{align*}
				\left\lvert 2a\alpha\left((b\beta + c\gamma)x
				+ (c\beta - b\gamma)y\right) \right\rvert &\le 2\lvert a\alpha \rvert \lvert b\beta x+c\beta y\rvert + 2\lvert a\alpha \rvert \lvert c\gamma x -b\gamma y\rvert \\
				&\le a^{2}\beta^{2} + \alpha^{2}(bx+cy)^{2} + a^{2}\gamma^{2} +\alpha^{2}(cx-by)^{2}\\
				&\le a^{2}(\beta^{2}+\gamma^{2}) +\alpha^{2}(b^{2}+c^{2}).
			\end{align*}
		\end{proof}

			This allows us to conclude that $R(Y \wedge JY,\Xi \wedge J\Xi)$ is nonpositive. Easy case by case examinations show that it vanishes only if $Y=0$ or $\Xi=0$. In conclusion, $\mu_{f,g}$ has negative holomorphic bisectional curvature on the cusp, and thus nonpositive holomorphic bisectional curvature on the compactified cusp. Gluing these metrics to the complex hyperbolic one in the thick part of $\HNC/\Gamma$ gives a Kähler metric on $X_{\Gamma}$ of nonpositive holomorphic bisectional curvature. This concludes the proof of Theorem \ref{thm_2_metrique_HBC}.
\

\begin{remark}
	In fact the metric constructed in the proof has negative Ricci curvature. The Ricci curvature, being a sum of holomorphic bisectional curvatures, is negative on the open subset $\HNC/\Gamma \subset X_{\Gamma}$, and we will prove that it is negatively pinched near the boundary tori. More precisely, using the notation of the proof of Theorem \ref{thm_2_metrique_HBC}, we claim that $\Ricci \le -2\mu_{f,g}$ on\linebreak $(0,\epsilon)\times N/\phi(\Lambda)$ for some $\epsilon>0$. Indeed, let $(X_1,\dots,X_{n-1})$ be an orthogonal $\C$-basis of $(\mathfrak{r},\mu\vert_{ \mathfrak{r}})$, so that, implictely evaluating $f$ and $g=ff'$ at some real number $t\in (0,A]$, the family\begin{equation*}
		\left(\frac{1}{f}X_1,\frac{1}{f}JX_1,\dots,\frac{1}{f}X_{n-1},\frac{1}{f}JX_{n-1},\frac{1}{g}Z,\frac{1}{g}JZ\right)
	\end{equation*}
	forms an orthonormal basis of the tangent space of $(0,A]\times N/\phi(\Lambda)$ at $(t,\Id \mod \varphi(\Lambda))$. Let
	
	\begin{equation*}
		\Xi = \sum_{i=1}^{n-1}(\alpha_iX_i+\alpha_{\bar{i}}JX_i) + \beta Z + \gamma JZ\\
	\end{equation*}
	be a tangent vector at this point, with $\alpha_i,\alpha_{\bar{i}},\beta,\gamma \in \R$, and denote by $\alpha$ the real number defined by\begin{equation*}
		\alpha := \sqrt{\sum_{i=1}^{n-1}(\alpha_i^{2}+\alpha_{\bar{i}}^{2})}.
	\end{equation*}
	Using the above formula for the bisectional curvature, a computation shows that the Ricci curvature of the metric has the following expression:
	\begin{equation*}
		\Ricci(\Xi,\Xi)=-\left(2ff''+4f'^{2}+2(n-2){f'^{2}}\right)\alpha^{2} -
		\left((2n+1)ff'^{2}f''+f^{2}f'f'''\right)(\beta^{2}+\gamma^{2}).
	\end{equation*}
	The function $f$ coincides with $\cosh$ on $(0,\epsilon)$, for some $\epsilon>0$ and when $t<\epsilon$, we get
	\begin{equation*}
		\Ricci(\Xi,\Xi) \le -2f^{2}\alpha^{2}-2g^{2}(\beta^{2}+\gamma^{2}) = -2\lVert \Xi\rVert^{2}_{\mu_{f,g}}.
	\end{equation*}
	By continuity, this inequality remains valid on the boundary tori, hence the result.
\end{remark}

		\section{Commensurators, parabolic subgroups and Albanese maps}\label{section_imm_alb}
		
		We begin this section by giving a definition and stating a proposition which will be used during the proof of Theorem \ref{thm_3_imm_alb}. 
		
		\begin{definition}
			For $\xi\in \PHNC$ and $o\in\HNC$, we define $K(\xi,o)$ as the subgroup of all elements in $\PU(n,1)$ fixing both $\xi$ and $o$.
		\end{definition}
		
		We first recall some facts about the group $K(\xi,o)$ that will be needed later. First this group is isomorphic to $\mathrm{U}(n-1)$. Moreover, recall from  Section \ref{section_preli} that $N_{\xi}$ is the unipotent radical of the stabilizer of $\xi$ in $\PU(n,1)$. As explained there, if $Z_{\xi}$ is a non-trivial element in the centre of $N_{\xi}$, then $\HNC/\langle Z_{\xi}\rangle$ is biholomorphic to the open subset $\Omega$ of $\C^{*}\times \C^{n-1}$ defined by Formula \eqref{eq-def-omega}. Define also\begin{equation*}
			\widehat{\Omega} := \Omega\cup (\{0\}\times \C^{n-1}).
		\end{equation*}
		It can be verified that the action of $K(\xi,o)$ on $\HNC$ commutes with $\langle Z_{\xi}\rangle$ so $K(\xi,o)$ also acts on $\Omega$. Identifying $K(\xi,o)$ with $\mathrm{U}(n-1)$, this action can simply be written $A\cdot(\lambda,v)=(\lambda,Av)$ for all $(\lambda,v) \in \Omega \subset \C^{*}\times \C^{n-1}$. This follows from Formula \eqref{eq_coor_Ltilde} and Diagram \eqref{diagramme}, and from the fact that $K(\xi,o)$ acts only on the ``$v$-terms'' of the coordinates $\widetilde{\mathcal{L}}$ and preserves its norm. In particular, the action of $K(\xi,o)$ extends by the same formula to $\widehat{\Omega}$.\newline
		
		We now turn to Proposition \ref{prop_comm_dense}. This proposition is probably classical, and we include it for the sake of being self-contained. If $\Gamma$ is a discrete subgroup of $\PU(n,1)$, recall that its commensurator $\Comm(\Gamma)$ is the set of elements $g\in \PU(n,1)$ such that $g\Gamma g^{-1}\cap \Gamma$ has finite index both in $\Gamma$ and in $g\Gamma g^{-1}$; it is a subgroup of $\PU(n,1)$. Moreover a parabolic point of $\Gamma$ is a point $\xi \in \PHNC$ whose stabiliser $\Stab_{\Gamma}(\xi)$ in $\Gamma$ contains a parabolic element.
		
		\begin{proposition}\label{prop_comm_dense}
			Let $\Gamma_0$ be a non-uniform arithmetic lattice in $\PU(n,1)$. Then there exists a dense subset $S\subset\HNC$ such that for any parabolic point $\xi$ of $\Gamma_0$ and for any point $o\in S$, $\Comm(\Gamma_0)\cap K(\xi,o)$ is dense in $K(\xi,o)$.
		\end{proposition}
		
		We defer the proof to the end of this section, and continue with Theorem \ref{thm_3_imm_alb}. For every closed Kähler manifold $X$, one can construct its Albanese torus $\Alb(X)$ and a holomorphic map from $X$ to $\Alb(X)$, called the Albanese map of $X$ \cite[chapter 12]{voisinTheorieHodgeGeometrie2002}. We denote by $\Alb(X_\Gamma)$ the Albanese torus of $X_{\Gamma}$ and $A_{\Gamma}:X_\Gamma \to \Alb(X_{\Gamma})$ its Albanese map. With this notation, let us restate Theorem \ref{thm_3_imm_alb}.
		
		\thmTrois*
		
			We recall that  by ``boundary torus of $X_{\Gamma}$'', we mean one of the tori added during the compactification, i.e. a connected component of the complement of $\HNC/\Gamma$ in $X_{\Gamma}$. We also mention that in dimension $n\ge3$, the boundary tori of the toroidal compactification $X_{\Gamma}$ are never simple when the lattice $\Gamma$ is arithmetic: see the remark at the end of this section.
		
		\begin{proof}
			For the first point, we refer to \cite[Theorem 1.3]{eyssidieuxOrbifoldKahlerGroups2018}, which uses that there are non-trivial holomorphic one forms on $X_{\Gamma}$ by results of Shimura \cite{shimuraAutomorphicFormsPeriods1979-manuel}. See also \cite[Section 3.1]{llosaisenrichSubgroupsHyperbolicGroups2024} for another exposition of this proof in the cocompact case.
			For the proof of the second point, we can and do assume that $\HNC/\Gamma_0$ admits a toroidal compactification, and we introduce the following terminology: a boundary torus $\tau$ of $X_{\Gamma_0}$ will be said to be \textit{virtually Albanese non-constant} if there exists a finite index subgroup $\Gamma < \Gamma_0$ and a boundary torus $T$ of $X_{\Gamma}$ above $\tau$, in restriction to which $A_{\Gamma}$ is not constant. 	
			Let $\tau$ be a virtually Albanese non-constant boundary torus of $X_{\Gamma_0}$, and $\Gamma$ and $T$ be as above. Holomorphic maps between tori have constant rank, thus for every point $x\in T$, there exists a holomorphic one form $\alpha \in \Omega^{1}(X_{\Gamma})$ such that $i^{*}\alpha(x)\ne 0$, where $i:T\hookrightarrow X_{\Gamma}$ is the inclusion.\newline
			
			Let $\xi\in \PHNC$ be a point corresponding to the cusp of $\HNC/\Gamma$ compactified by the torus $T$. From here on, we use the notation of Section \ref{section_preli}. According to the results recalled there, $\Stab_{\Gamma}(\xi)$ is a lattice in the nilpotent group $N_{\xi}$. Thus it intersects $[N_{\xi},N_{\xi}] \simeq \R$ non-trivially. Let $Z_{\xi}$ be a generator of $\Gamma\cap [N_{\xi},N_{\xi}]$. We have the following diagram
			
			\begin{equation*}
				\begin{tikzcd}
					\HNC/\langle Z_{\xi}\rangle \arrow{d}{} \arrow{r}{\mathcal{L}}  & \Omega\arrow[hookrightarrow]{r} &\widehat{\Omega} \arrow{d}\arrow[hookleftarrow]{r}{\widetilde{i}} &\C^{n-1}\arrow{d}\\
					\HNC/\Gamma \arrow[hookrightarrow]{rr} && X_{\Gamma} \arrow[hookleftarrow]{r}{i} &T,
				\end{tikzcd}
			\end{equation*}
			where $\widehat{\Omega} := \Omega\cup (\{0\}\times \C^{n-1})$, hooked arrows are inclusions, $\widetilde{i}$ is defined by $\widetilde{i}(x):=(0,x)$, and the map $\widehat{\Omega}\to X_{\Gamma}$ is a covering above $T\cup \HNC/\Gamma$. Recall that $\mathcal{L}$ depends on the choice of a geodesic $\gamma$ such that $\gamma(t)\underset{t\to+\infty}{\longrightarrow}\xi$. Here we choose $\gamma$ such that $o:=\gamma(0)$ belongs to the set $S$ of Proposition \ref{prop_comm_dense}. As explained in the beginning of this section, the group $K(\xi,o)$ acts on $\widehat{\Omega}$.\newline
			
			Let $x\in X_{\Gamma}$ be the point below $(0,0)\in \widehat{\Omega}$, and $\alpha \in \Omega^{1}(X_{\Gamma})$ a holomorphic one form such that $i^{*}\alpha(x)\ne 0$. Denoting by $\beta$ the pullback of $\alpha$ by the map $\widehat{\Omega}\to X_{\Gamma}$, we have $\widetilde{i}^{*}\beta(0)\ne 0$.	We deduce that there exist $g_1,\dots,g_{n-1} \in K(\xi,o)$ such that \begin{equation*}
				e: =(\widetilde{i}^{*}(g_{1}^{-1})^{*}\beta(0), \dots, \widetilde{i}^{*}(g_{n-1}^{-1})^{*}\beta(0))
			\end{equation*} is a basis of $T^{(1,0)*}_{0}\C^{n-1}$. All the $g_i$'s can be chosen in $\Comm(\Gamma)\cap K(\xi,o)$, which is dense in $K(\xi,o)$ by Proposition \ref{prop_comm_dense}.\newline
			
			Then let $\Gamma'$ be some finite index subgroup of\begin{equation*}
				\bigcap_{i=1}^{n-1}g_{i}\Gamma g_{i}^{-1} \cap \Gamma
			\end{equation*}
			 which is normal in $\Gamma$. We want to show that for every boundary torus $T'$ of $X_{\Gamma'}$ above $T$ and for every point $x'\in T'$ above $x$, the differential $dA_{\Gamma'\,\vert T'}(x')$ is injective. Since $\Gamma'$ is normal in $\Gamma$, the automorphism group of $X_{\Gamma'}$ acts transitively on the set of tori above $T$, thus it is enough to show this property for the torus $T'$ above $T$ which compactifies the cusp associated with the $\Gamma'$-orbit of $\xi$.\newline
			
			The group $\Gamma' \cap [N_{\xi},N_{\xi}]$ is a non-trivial subgroup of $\Gamma \cap [N_{\xi},N_{\xi}]$, it is thus generated by $dZ_{\xi}$ for some $d\in \N^{*}$. We then have the following diagram
			
			\begin{equation}\label{diag_comm_rev}
				\begin{tikzcd}
					\HNC/\langle dZ_{\xi}\rangle \arrow{d}{} \arrow{r}{\mathcal{L}^{(d)}}  & \Omega^{(d)}\arrow[hookrightarrow]{r} &\widehat{\Omega^{(d)}} \arrow{d}{\Psi}\arrow[hookleftarrow]{r}{\widetilde{i}^{(d)}} &\C^{n-1}\arrow{d}{\Id}\\
					\HNC/\langle Z_{\xi}\rangle \arrow{d}{} \arrow{r}{\mathcal{L}}  & \Omega\arrow[hookrightarrow]{r} &\widehat{\Omega} \arrow{d}\arrow[hookleftarrow]{r}{\widetilde{i}} &\C^{n-1}\arrow{d}\\
					\HNC/\Gamma \arrow[hookrightarrow]{rr} && X_{\Gamma} \arrow[hookleftarrow]{r}{i} &T,
				\end{tikzcd}
			\end{equation}
			where\begin{align*}
				&\Omega^{(d)} := \left\{(b,v) \in \C\times \C^{n-1} \ \vert\ 0<\lvert b\rvert < \exp\left({\frac{-\pi \lVert v\rVert^{2}}{dl}}\right)\right\},\\
				&\widehat{\Omega^{(d)}} := \Omega^{(d)}\cup (\{0\}\times \C^{n-1}) \qquad \text{ and }\\
				&\Psi\colon (b,v) \mapsto (b^{d},v).
			\end{align*}
			
			Let $\delta:=\Psi^{*}\beta$, and define $\delta_i := (g_i^{-1})^{*}\delta$ for $i\in\{1,\dots,n-1\}$. We claim that $\delta_1,\dots,\delta_{n-1}$ are the pullback of holomorphic forms $\epsilon_1,\dots,\epsilon_{n-1}$ on $X_{\Gamma'}$, and that denoting $T'\overset{i'}{\hookrightarrow} X_{\Gamma'}$ the inclusion of $T'$ in $X_{\Gamma'}$, the family $e'':=(i'^{*}\epsilon_1(x'),\dots,i'^{*}\epsilon_{n-1}(x'))$ is a basis of $T^{(1,0)*}_{x'}T'$. This will imply that the differential $dA_{\Gamma'\,\vert T'}(x')$ is injective.\newline
			
			For all $i$, the form $\delta_i$ is the pullback of a holomorphic form $\epsilon_i$ because the element $g_i \in \PU(n,1)$ induces a biholomorphism $\HNC/\Gamma \to \HNC/g_i\Gamma g_i^{-1}$ which extends to a biholomorphism $\theta_i\colon X_{\Gamma} \to X_{g_i\Gamma g_i^{-1}}$. Noticing that $g_i$ acts on $\widehat{\Omega^{(d)}}$ and that there is a ramified covering \linebreak $\pi_i:X_{\Gamma'}\to X_{g_i\Gamma g_i^{-1}}$, the commutativity of the following diagram entails that $\delta_i$ is the pullback of the form $\epsilon_i := \pi_i^{*}(\theta_i)_*\alpha$
			\begin{equation*}
				\begin{tikzcd}
					\widehat{\Omega^{(d)}} \arrow{d}{\Psi} \arrow{r}{g_i\cdot} &\widehat{\Omega^{(d)}} \arrow{r}&X_{\Gamma'} \arrow{dd}{\pi_i} \\
					\widehat{\Omega}\arrow{d} &&\\
					X_{\Gamma} \arrow{rr}{\theta_i}&&X_{g_i\Gamma g_i^{-1}}.
				\end{tikzcd}
			\end{equation*}
			The family $e''$ is a basis because, up to the cover $\C^{n-1}\to T'$ it can be identified with the family $e':=(\widetilde{i}^{(d)*}\delta_1(0),\dots,\widetilde{i}^{(d)*}\delta_{n-1}(0))$ which, up to the identity map $\C^{n-1}\to \C^{n-1}$ of the commutative diagram \eqref{diag_comm_rev}, is (identified with) $e$. The latter was constructed to be a basis of $T^{(1,0)*}_{0}\C^{n-1}$.\newline
			
			We deduce that $dA_{\Gamma' \,\vert T'}(x')$ is injective. Since $A_{\Gamma' \,\vert T'}$ is a holomorphic map between tori, it is a fibration, and in particular, its differential is injective at a point if and only if it is injective at any point. Hence $A_{\Gamma' \,\vert T'}$ is an immersion.\newline
			
			For each virtually Albanese non-constant boundary torus $\tau_i$, there is a finite index subgroup $\Gamma'_i<\Gamma_0$ constructed as above. Let $\Gamma'$ be the intersection of the $\Gamma'_i$. Let $\Gamma''$ be a finite index subgroup of $\Gamma'$, and $T$ a boundary torus of $X_{\Gamma''}$. If $T$ covers a boundary torus of $X_{\Gamma_0}$ which is not virtually Albanese non-constant, then by definition $A_{\Gamma'' \,\vert T}$ is constant. Otherwise, there exists an intermediate ramified covering $X_{\Gamma''} \to X_{\Gamma'_i}\to X_{\Gamma_0}$ such that the Albanese map of $X_{\Gamma'_i}$ in restriction to a torus below $T$ is an immersion. \textit{A fortiori} the Albanese map of $X_{\Gamma''}$ in restriction to $T$ is an immersion.
		\end{proof}
		
		\begin{proof}[Proof of Proposition \ref{prop_comm_dense}]
			We use here the explicit description of non-uniform arithmetic lattices in $\PU(n,1)$, see for instance \cite[Section 2]{emeryCovolumesNonuniformLattices2014} or \cite{prasadArithmeticFakeProjective2009-manuel}. We denote $\pi:\U(n,1)\to \PU(n,1)$ and $[\cdot]:\C^{n+1}\setminus\{0\}\to \PNC$ the projection maps. Let $\Gamma_0$ be a non-uniform arithmetic lattice in $\PU(n,1)$. Then there exist:\begin{itemize}
				\item A purely imaginary quadratic extension $l:=\mathbb{Q}(\sqrt{-d})$ of $\mathbb{Q}$ whose ring of integers is denoted $\O_l$,
				\item A Hermitian form $H:l^{n+1}\times l^{n+1}\to l$ on $l^{n+1}$ whose extension $H_{\C}$ to $\C^{n+1}$ is of signature $(n,1)$,
				\item An isomorphism $\Phi:(\C^{n+1},H_{\C})\to (\C^{n+1},\langle \cdot,\cdot\rangle)$ where $\langle \cdot,\cdot\rangle$ denotes the standard Hermitian form of signature $(n,1)$ used to define $\HNC$; moreover $\Phi$ induces an isomorphism $\Theta: \U(H_{C})\to \U(n,1)$ defined by $\Theta(g):=\Phi g\Phi^{-1}$, and allows to define an arithmetic lattice $\Gamma_l<\PU(n,1)$ by $\Gamma_l := \pi(\Theta(U(H_{C},\O_l)))$; and
				\item A subgroup $\Gamma$ of $\Gamma_0\cap\Gamma_l$ which is of finite index both in $\Gamma_0$ and $\Gamma_l$.
			\end{itemize}
			If we replace $\Gamma$ by a finite index subgroup, we can also suppose that all its parabolic subgroups are unipotent. We define\begin{align*}
				\HNL := \{[\Phi(x)] \ \vert \ x\in l^{n+1} \text{ and } H(x,x)<0\},\\
				\PHNL := \{[\Phi(x)] \ \vert \ x\in l^{n+1} \text{ and } H(x,x)=0\}.
			\end{align*}
			Let $\xi:=[\Phi(\widetilde{\xi})] \in \PHNL$ and $o:=[\Phi(\widetilde{o})]\in \HNL$, with $(\widetilde{\xi},\widetilde{o})\in l^{n+1}\times l^{n+1}$. Let $C$ be the intersection of the stabilisers of $\widetilde{\xi}$ and $\widetilde{o}$ in $U(H_{\C})$, and $C_{l}$ the subgroup of $C$ consisting of the matrices with coefficients in $l$. By construction, $K(\xi,o)=\pi(\Theta(C))$.
			
			Two commensurable lattices have the same commensurator and the same parabolic points. Therefore, the following three statements lead to the proposition with $S:=\HNL$.
			
			\underline{Assertion 1}: let $\xi$ be a parabolic point of $\Gamma$. Then $\xi\in \PHNL$.
				\vspace*{.2\baselineskip}
				
			\underline{Assertion 2}: let $\xi\in\PHNL$ and $o\in \HNL$. Then $\pi(\Theta(C_{l})) \subset \Comm(\Gamma_l)\cap K(\xi,o)$, with $C_l$ defined as above.
				\vspace*{.2\baselineskip}
				
			\underline{Assertion 3}: let $\xi\in \PHNL$ and $o\in \HNL$. Then $C_{l}$ is dense in $C$, so $\pi(\Theta(C_{l}))$ is dense in $K(\xi,o)$.
			
			The proofs of these three statements are independent of each other. All these assertions are very classical. We give elementary proofs of Assertions 1 and 3, though shorter and more general proofs exist.\newline
			
			\textit{Proof of Assertion 1 (see also \cite{zinkUeberAnzahlSpitzen1979})}.
			
			Let $\gamma$ be a non-trivial parabolic element of $\Gamma$ which fixes $\xi=[x]$. We must show that $\C x$ intersects non trivially $l^{n+1}$. Since $\Gamma\subset \Gamma_l$, we can write $\gamma=\pi(\Theta(M))$ with $M\in \U(H_{\C},\O_l)$, and $\Theta(M) x=\alpha x$ for some $\alpha\in S^{1}$. Since $\gamma$ is unipotent, it can also be represented by an upper triangular matrix $T$ with diagonal entries 1, thus $\Theta(M)$ is a multiple $\alpha'T$ of $T$. Evaluating at $x$ gives $\alpha=\alpha'$, and in particular $\det(M)=\det(\Theta(M))=\alpha^{n+1}$. We deduce that $\alpha^{n+1}\in \O_l$. The only integers of modulus 1 of purely imaginary quadratic extensions are roots of units, so $\alpha$ is a root of unit. Up to replacing $\gamma$ by an appropriate power $\gamma^{k}$, we can thus suppose that $\gamma=[\Theta(M)]$ with $\Theta(M)x=x$.
			
			The map $M$ can also be seen as a linear map $M_l:l^{n+1}\to l^{n+1}$. The kernel $E^{1}_l$ of $M_l-\Id$ and the kernel $E^{1}$ of $M-\Id$ are related by $E^{1}=E^{1}_l\otimes \C$. Let $e:=(e_1,\dots,e_k)$ be an orthogonal basis of $E^{1}_l$. It is not possible to have $H(e_i,e_i)>0$ for all $i$ because otherwise, $H_{\C}$ would be positive definite on the $\C$-span of $e$, which is $E^{1}$, and this space contains the isotropic vector $x$. Thus $H(e_i,e_i)\le 0$ for some $i$. Now $\xi=[x]$ and $[\Phi(e_i)]$ are two fixed points  in $\overline{\HNC}$ of the parabolic element $\gamma=\pi(\Theta(M))$, hence $\xi=[\Phi(e_i)]\in \PHNL$.\newline

			\textit{Proof of Assertion 3}.
			
			Assertion 3 is a consequence of the fact that $C_l$, being an algebraic group defined over $\mathbb{Q}$, is unirational over $\mathbb{Q}$ \cite[Theorem 18.2]{borelLinearAlgebraicGroups1991}. We now give a self-contained proof of this assertion.
			
			Let $\widetilde{\xi}\in l^{n+1}$ and $\widetilde{o}\in l^{n+1}$ be respectively isotropic and negative vectors of $H$, $P$ the $l$-span of these two vectors and $P^{\perp}$ its orthogonal with respect to $H$. Then, since $H$ has signature $(n,1)$, $H(\widetilde{\xi},\widetilde{o})\ne 0$ and this easily implies that $P\cap P^{\perp}=\{0\}$. Hence $l^{n+1}=P\oplus P^{\perp}$, and $H$ restricted to $P^{\perp}$ is a positive definite Hermitian form. In particular, there is a basis $(u_1,\dots,u_{n-1})$ of $P^{\perp}$ orthogonal for $H$, and the matrix $B:=(H(u_i,u_j))_{i,j}$ is diagonal, with positive rational diagonal entries. Now define\begin{equation*}
				\U := \{M \in \Mat(n-1,\C) \ \vert\ ^{t}MB\overline{M} = B\}.
			\end{equation*}
			Then $C$ is isomorphic to $\U$ by a Lie group isomorphism sending $C_l$ to\begin{equation*}
				\U_l := \{M \in \Mat(n-1,l) \ \vert\ ^{t}MB\overline{M} = B\}.
			\end{equation*}
			We will show that $\U_l$ is dense in $\U$, using a classical argument involving the Cayley transform \cite[section II.10]{weylClassicalGroupsTheir1946}. Here are the notation that we are going to use:\begin{align*}
					&\Mat^{*} := \{M\in \Mat(n-1,\C) \ \vert\ \det(M+\mathrm{I})\ne 0\},\\
					&\Mat_l^{*} := \Mat^{*} \cap \Mat(n-1,l),\\
					&\U^{*} := \U\cap\Mat^{*},\\
					&U^{*}_l := \U^{*}\cap \Mat(n-1,l),\\
					&\mathrm{A}^{*} := \{M \in \Mat(n-1,\C) \ \vert\  ^{t}MB=-B\overline{M}\} \cap \Mat^{*}, \\
					& \mathrm{A}^{*}_l := \mathrm{A}^{*}\cap \Mat(n-1,l).
			\end{align*}
			The Cayley transform is the involution\begin{align*}
				\begin{array}{rrcl}
					S\colon &\Mat^{*} &\longrightarrow &\Mat^{*}\\
					& N &\longmapsto & 2(\mathrm{I}+N)^{-1}-\mathrm{I}.
				\end{array}
			\end{align*}
			It is easily verified that $S:\U^{*}\to \mathrm{A}^{*}$ is a homeomorphism sending $\U^{*}_l$ to $\mathrm{A}^{*}_l$. Since $\U^{*}$ is dense in $\U$, it is thus enough to show that $\mathrm{A}^{*}_l$ is dense in $\mathrm{A}^{*}$. This follows from the fact that $A^{*}$ is the intersection of a dense open set and a subspace defined by rational equations.
		\end{proof}
		
		\begin{remark}
			We outline there a proof of the classical fact that the boundary tori arising from arithmetic toroidal compactifications in dimension $n\ge 3$ are not simple. We will use the description of arithmetic lattices given in the proof of Proposition \ref{prop_comm_dense}. Let $\Gamma_0$ be an arithmetic lattice in $\PU(n,1)$, $\Gamma$ be a finite index subgroup of $\Gamma_{0}$ as in the proof of Proposition \ref{prop_comm_dense} and $T$ be a boundary torus of $X_{\Gamma}$ which corresponds to a parabolic point $\xi\in\PHNC$. There is a basis $f_{\xi}=(f_1,f_2,\dots,f_{n+1})$ of $\C^{n+1}$ as in Section \ref{section_preli} with coefficients in $\ell$. In this basis, the matrix group that corresponds to the parabolic subgroup of $\Gamma$ fixing $\xi$ is of the form $\widehat{P}:=gPg^{-1}$ with $g\in \mathrm{GL}_{n+1}(\ell)$ and $P$ a subgroup of $\mathrm{GL}_{n+1}(\mathcal{O}_\ell)$. We deduce that there is a finite-index subgroup $\widehat{P}_0$ of $\widehat{P}$ that is contained in $\mathrm{GL}_{n+1}(\mathcal{O}_\ell)$. The quotient of $\widehat{P}_0$ by its center has a a finite index subgroup of the form $H^{n-1}$, for some finite index subgroup $H$ of $\mathcal{O}_\ell$. In other words, $T$ has a finite cover that is a product of elliptic curves, and we conclude that $T$ is not simple.
		\end{remark}
		
		\section{The Shafarevich conjecture for the manifold $X_\Gamma$}\label{section_rev_stein}
		
		In this section we prove Theorem \ref{thm_4_rev_Stein}. We describe the global strategy of the proof in the beginning of Subsection \ref{ss-section-univ_Stein}. The results which apply in a more general context and may be of independent interest have been gathered in Subsection \ref{sssection_res_gen}. In Subsection \ref{ss-section-univ_Stein}, we return to toroidal compactifications and show that the Albanese map of $X_\Gamma$ satisfy the hypotheses of Proposition \ref{lem-PSH3} stated and proven in Subsection \ref{sssection_res_gen}.
		
		\subsection{A general result on holomorphic maps whose target have Stein universal cover}\label{sssection_res_gen}
		
		 The following proposition is stated for projective manifolds which are the domain of a holomorphic map with values in a compact manifold with Stein universal cover. Lemma \ref{lem-PSH4} is used to prove the second point of the proposition. We first give a definition which will be used through all Section \ref{section_rev_stein}.
		 
		 \begin{definition}
		 	Let $V$ be a subset of a topological space $X$ and $f:V\to \R$ be a continuous function. Then $f$ is said to be an \textit{exhaustion relative to }$X$ if for all real number $\lambda$, the closure in $X$ of the sublevel set $f^{-1}((-\infty,\lambda))$ is compact. 
		 \end{definition}
		
		\begin{proposition}\label{lem-PSH3}
			Let $X$ be a projective manifold, $\pi:\widetilde{X}\to X$ its universal cover and $A:X\to Y$ a holomorphic map from $X$ to a complex manifold $Y$ whose universal cover $\widetilde{Y}$ is Stein. \begin{enumerate}
				\item Define\begin{equation*}
					Z := \{x\in X \ \vert \ A^{-1}(A(x)) \text{ is not finite}\},
				\end{equation*}
				let $U$ be an open neighborhood of $Z$ and set $\widetilde{U}:=\pi^{-1}(U)$. Then $Z$ is an analytic subvariety of $X$ and there is a smooth plurisubharmonic function $\psi\colon \widetilde{X} \to \R_+$, whose vanishing locus is exactly $\widetilde{Z}:=\pi^{-1}(Z)$, that is strictly plurisubharmonic on $\widetilde{X}\setminus \widetilde{Z}$, and whose restriction to $\widetilde{X}\setminus \widetilde{U}$ is an exhaustion.
				\item In addition, if there is an open neighborhood $V$ of $Z$ and a smooth strictly plurisubharmonic function $\phi$ on $\widetilde{V}:=\pi^{-1}(V)$ which is an exhaustion relative to $\widetilde{X}$, then $\widetilde{X}$ is Stein.
			\end{enumerate}
		\end{proposition}

		\begin{lemma}\label{lem-PSH4}
			Let $X$ be a metric space and $U,V$ two open subsets of $X$ with $\overline{U}\subset V$. Suppose that there exist:\begin{itemize}
				\item A continuous function $\phi:V\to (0,+\infty)$, and
				\item A continuous function $\psi:X\to \R_+$ such that $\psi^{-1}(\{0\})\subset U$ and $\psi\vert_{ X\setminus U}$ is an exhaustion.
			\end{itemize}
			Then for any open subset $V'$ containing $U$ whose closure is included in $V$, there exists an increasing smooth convex function $\chi:\R_+\to \R_+$ such that $\chi(0)=0$ and $\chi\circ\psi>\phi$ on $V'\setminus U$.
		\end{lemma}
		
		\begin{proof}[Proof of Proposition \ref{lem-PSH3}-1] This proof is mainly taken from \cite[Proof of Theorem 6.1]{napierConvexityPropertiesCoverings1990}, with some modifications to take into account that we work in arbitrary dimension.
			
			\underline{Step 1.} The set $Z$ is an analytic subvariety of $X$.
			
			This is classical, see for instance \cite[§3.6]{fischerComplexAnalyticGeometry1976}. For a more elementary argument using that $X$ is projective, one can also realise $X$ as a submanifold of some projective space $\C\mathbb{P}^{d}$ and then choose for all $x \in X\setminus Z$ a hyperplane $H_x$ of $\C\mathbb{P}^{d}$ which does not intersect the finite set $A^{-1}(A(x))$. Then $Z$ is analytic as a consequence of the following equality, which is easily verified:\begin{equation*}
				Z = \bigcap_{x \in X\setminus Z}A^{-1}(A(H_x)).
			\end{equation*}
		\
		
			\underline{Step 2.} For any sequence $(\widetilde{x_k})_k$ of $\widetilde{X}\setminus\widetilde{U}$ with no accumulation point, there is a holomorphic function $f:\widetilde{X} \to \C$ which is not bounded on $(\widetilde{x_k})_k$ and vanishes on $\widetilde{Z}:=\pi^{-1}(Z)$.
			
			Let $\rho:\widetilde{Y}\to Y$ be the universal cover of $Y$ and $\widetilde{A}\colon \widetilde{X} \to \widetilde{Y}$ the lift of $A$, so that the following diagram commutes:
			\begin{equation*}
				\begin{tikzcd}
					\widetilde{X}\arrow{d}{\pi}\arrow{r}{\widetilde{A}}&\widetilde{Y}\arrow{d}{\rho} \\
					X \arrow{r}{A}  & Y.
				\end{tikzcd}
			\end{equation*}
			If $(\widetilde{A}(\widetilde{x_k}))_k$ has no accumulation point, then as $\widetilde{Y}$ is a Stein manifold, there is a holomorphic function $u:\widetilde{Y}\to \C$ which is unbounded on $(\widetilde{A}(\widetilde{x_k}))_k$ and vanishes on the subvariety $\rho^{-1}(A(Z))$. The function $u\circ\widetilde{A}$ has the required properties.\newline
			
			Otherwise, the sequence $\widetilde{y_k} :=\widetilde{A}(\widetilde{x_k})$ has an accumulation point, and up to extracting a subsequence, we can assume, on the one hand, that it converges to an element $\widetilde{y_{\infty}} \in \widetilde{Y}$, and on the other hand that the sequence $x_k:=\pi(\widetilde{x_k})$ converges to a limit $x_{\infty} \in X\setminus Z$. We claim that there exists an ample line bundle $L\to X$ such that:\begin{itemize}
				\item There exists a holomorphic section $t\in H^{0}(\widetilde{X},\pi^{*}L)$ which is unbounded on $(\widetilde{x_k})_k$ with respect to $\pi^{*}h$ for some Hermitian metric $h$ on $L$.
				\item There exists a holomorphic section $s\in H^{0}(X,L)$ whose vanishing locus $S\subset X$ satisfies:
				\begin{itemize}
					\item $S\cap A^{-1}(A(x_{\infty}))=\emptyset$, in other words $A(x_{\infty}) \notin A(S)$. In particular, we have \linebreak $\pi^{*}s(\widetilde{x_k}) \longrightarrow s(x_{\infty})\neq 0$.
					\item The complement of $S$ does not contain a compact subvariety of positive dimension of $X$.
				\end{itemize}
			\end{itemize}
			Indeed, since $X$ is a projective manifold, it can be realized as a submanifold of some projective space $\C \mathbb{P}^{d}$. Let $H\subset \C \mathbb{P}^{d}$ be a projective hyperplane disjoint from the finite set $A^{-1}(A(x_{\infty}))$. Then $\mathcal{O}(1)\vert_{ X}$ is positive so by \cite[Corollary 4.3]{napierConvexityPropertiesCoverings1990}, for a large enough integer $p$, the bundle $L: =\mathcal{O}(p)\vert_{ X}$ is sufficiently positive so that there exists a holomorphic section $t\in H^{0}(\widetilde{X},\pi^{*}L)$ unbounded on $(\widetilde{x_k})_k$. For the second point, let $\widetilde{s_0}$ be a global holomorphic section of $\mathcal{O}(1)$ whose vanishing locus is $H$. Then the vanishing locus of the section $\widetilde{s}:=\widetilde{s_0}^{\otimes p}$ of $\mathcal{O}(p)$ is $H$ and its complement, biholomorphic to $\C^{d}$, is a Stein manifold and therefore contains no compact subvariety of positive dimension. The restriction $s \in H^{0}(X,L)$ of $\widetilde{s}$ to $X$, whose vanishing locus is $S:=X\cap H$ has the required properties. Set also $\widetilde{S}:=\pi^{-1}(S)$.\newline
			
			The set $E:=\rho^{-1}(A(S))$ is an analytic subset of $\widetilde{Y}$, which is equal to $\rho^{-1}(\rho\circ\widetilde{A}(\widetilde{S}))$. Since $A(x_{\infty}) \notin A(S)$, the point $\widetilde{y_{\infty}}$, whose image by $\rho$ is $A(x_{\infty})$, does not belong to $E$. We deduce that there exists a holomorphic function $u: \widetilde{Y}\to \C$ which vanishes on $E$ but does not vanish on $\widetilde{y_{\infty}}$. Indeed the sheaf $\mathcal{F}$ of holomorphic functions vanishing on $E$ is coherent, so by Cartan's Theorem B, there exist holomorphic functions vanishing on $E$ which generate the germ space $\mathcal{F}_{\widetilde{y_{\infty}}}=\mathcal{O}_{\widetilde{y_{\infty}}}$. Then, $E$ contains $\widetilde{A}(\widetilde{S})$, and, we will show that it also contains $\widetilde{A}(\widetilde{Z})$. Let $\widetilde{x}\in \widetilde{Z}$ and $x:=\pi(\widetilde{x})$. By definition of $Z$, the analytic subset $A^{-1}(A(x))$ is not finite, so has a non-empty intersection with $S$. In other words, $\rho\circ\widetilde{A}(\widetilde{x})=A(x) \in A(S)$, so $\widetilde{A}(\widetilde{x}) \in E$. Thus, $u\circ\widetilde{A}$ vanishes on $\widetilde{S}\cup\widetilde{Z}$. As the divisor of $\pi^{*}s$ is $p\widetilde{S}$, we deduce that the meromorphic function\begin{equation*}
				f:=(u\circ \widetilde{A})^{p+1}\times (t/\pi^{*}s)
			\end{equation*} 
			is in fact holomorphic. Since $t$ is not bounded on $(\widetilde{x_k})_k$, and  the sequences $u(\widetilde{y_k})$ and $\pi^{*}s(\widetilde{x_k})$ converge respectively to non-zero limits $u(\widetilde{y_{\infty}})$ and $s(x_{\infty})$, the function $f$ is not bounded on $(\widetilde{x_k})_k$. Finally it is divisible as a holomorphic function by $u\circ\widetilde{A}$, thus it vanishes on $\widetilde{Z}$.\newline
			
			\underline{Step 3.} For any point $\widetilde{x_0}\in \widetilde{X}\setminus\widetilde{Z}$, there exists a holomorphic function $F:\widetilde{X}\to \C^{n}$ which vanishes on $\widetilde{Z}$, does not vanish at $\widetilde{x_0}$ and whose differential at $\widetilde{x_0}$ is invertible ($n=\dim X$).
			
			Let us write $x_0:=\pi(\widetilde{x_0})$. We assert that there exists an ample line bundle $L\to X$ and holomorphic sections $s,t_1,\dots,t_n$ of $L$ such that:\begin{itemize}
				\item The section $s$ does not vanish on $A^{-1}(A(x_{0}))$, and in particular does not vanish at $x_0$ so that the meromorphic functions $t_i/s$ are holomorphic on a neighbourhood of $x_0$.
				\item In a neighbourhood of $x_0$, the holomorphic functions $(\frac{t_1}{s},\dots,\frac{t_n}{s})$ do not vanish and define coordinates of $X$. 
				\item The complement of the vanishing locus $S\subset X$ of $s$ does not contain a compact subvariety of positive dimension of $X$.
			\end{itemize}
			Indeed, embedding $X$ in $\C \mathbb{P}^{d}$ we choose $L:=\mathcal{O}(1)\vert_{ X}$, and $s:=\widetilde{s}\vert_{ X}$ where $\widetilde{s}$ is a section of $\mathcal{O}(1)$ whose vanishing hyperplane $H$ is disjoint from $A^{-1}(A(x_{0}))$. On the complement of $H$, the section $\widetilde{s}$ trivializes $\mathcal{O}(1)$, therefore it makes sense to speak of the differential of a section $\widetilde{t}$ of $\mathcal{O}(1)$ at a point of the complement of $H$: it is the differential of the function $\widetilde{t}/\widetilde{s}$. There are homogeneous coordinates $[z_1:\dots:z_{d+1}]$ for which $\widetilde{s}$ identifies with the linear form $z_{d+1}$, the point $x_0$ has coordinates $[1:1:\dots:1]$ and, in the affine coordinates given by \begin{equation*}
				(z_1,\dots,z_d) \mapsto [z_1:\dots:z_d:1],
			\end{equation*}the tangent space of $X$ at $x_0$ is spanned by $(\frac{\partial}{\partial z_1},\dots,\frac{\partial}{\partial z_n})$. The linear forms $z_1,\dots,z_n$ identify with global sections $\widetilde{t_1},\dots,\widetilde{t_n}$ of $\mathcal{O}(1)$. Let $t_1,\dots,t_n$ denote the restriction of these sections to $X$. By construction, the holomorphic functions $(\frac{t_1}{s},\dots,\frac{t_n}{s})$ do not vanish and define coordinates of $X$.
			
			Let $S:= X \cap H$ be the vanishing locus of $s$, and $\widetilde{S}:=\pi^{-1}(S)$. By construction, the functions $\pi^{*}(t_i/s)$ are meromorphic on $\widetilde{X}$, and their divisors satisfy $\mathrm{div}(\pi^{*}(t_i/s))\ge-\widetilde{S}$. As above, $E:=\rho^{-1} (\rho\circ\widetilde{A} (\widetilde{S}))$ is an analytic subset of $\widetilde{Y}$ containing $\widetilde{A}(\widetilde{S})\cup \widetilde{A}(\widetilde{Z})$ and not containing $\widetilde{A}(\widetilde{x_0})$. By Cartan's Theorem B, there exists a holomorphic function $u:\widetilde{Y}\to \C$ which vanishes on $E: = \rho^{-1} (\rho\circ\widetilde{A} (\widetilde{S}))$, does not vanish at $\widetilde{A}(\widetilde{x_0})$ and whose differential at $\widetilde{A}(\widetilde{x_0})$ is zero. In particular, $u\circ \widetilde{A}$ vanishes on $\widetilde{S}\cup \widetilde{Z}$ but not at $\widetilde{x_0}$. Let us set\begin{equation*}
				f_i := (u\circ \widetilde{A})^{2}\times \pi^{*}\frac{t_i}{s}.
			\end{equation*}
			Then $F:=(f_1,\dots,f_n)$ is holomorphic. We assert that $F$ vanishes on $\widetilde{Z}$, does not vanish at $\widetilde{x_0}$ and that its differential is invertible at $\widetilde{x_0}$. Indeed, for all $i$, $f_i$ is divisible by $u\circ \widetilde{A}$, thus it vanishes on $\widetilde{Z}$. However $f_i$ does not vanish at $\widetilde{x_0}$ by construction. Finally, the differential of $u\circ \widetilde{A}$ is zero at $\widetilde{x_0}$ so\begin{equation*}
				(df_i)_{\widetilde{x_0}} = u(\widetilde{A}(\widetilde{x_0}))^{2} \pi^{*}d\left(\frac{t_i}{s}\right)_{x_0}.
			\end{equation*}
			Therefore, the differential of $F$ is invertible at $\widetilde{x_0}$.
			
			\underline{Step 4.} Construction of $\psi$ from these two statements: we refer to \cite[pp. 470-472]{napierConvexityPropertiesCoverings1990}.
		\end{proof}
		
		To prove the second point of Proposition \ref{lem-PSH3}, we need Lemma \ref{lem-PSH4}, that we prove here.
		
		\begin{proof}[Proof of Lemma \ref{lem-PSH4}]
			First, we claim that $m:=\inf_{V'\setminus U}(\psi)$ is positive. Indeed, seeking a contradiction, suppose that $m=0$. Then there is a sequence $(v_n)$ of elements of $V'\setminus U$ such that $\psi(v_n)\to 0$. By assumption, $\psi\vert_{ X\setminus U}$ is an exhaustion, thus up to taking some subsequence, we can assume that $(v_n)$ converges toward a limit $v_{\infty}\in X\setminus U$ such that $\psi(v_{\infty})=0$. This contradicts the hypothesis $\psi^{-1}(\{0\})\subset U$.
			
			If $\phi$ is bounded on $V'\setminus U$ by some constant $C$, then $\chi(x):=\frac{C+1}{m}x$ is suitable. Assume now that $\phi$ is unbounded on $V'\setminus U$, and for all $p\in \N^{*}$, define the set $A_p := \{x\in V'\setminus U \ \vert\ p-1\le \phi(x)< p\}$. Enumerate the integers $p$ such that $A_{p}\ne \emptyset$ as an increasing sequence of integers $(p_k)_{k\in \N}$. For all integer $k$, let $a_{k}$ be the infimum of the function $\psi\vert_{ \overline{A_{p_k}}}$, and $x_{k}\in \overline{V'}\setminus U$ be a point such that $\psi(x_k)\le a_k+1$.
			We claim that $a_k\underset{k\to +\infty}{\longrightarrow}+\infty$. Indeed, seeking a contradiction, suppose that there is a subsequence $(a_{k_j})_{j\in \N^{*}}$ of $(a_k)$, bounded by some constant $M$. Then the sublevel $\{\psi\vert_{ X\setminus U} \le M+1\}$ is compact and contains $x_{k_j}$ for all $j$ so up to a subsequence, we can assume that $x_{k_j}$ converges to a limit $\ell \in \overline{V'}$. We obtain a contradiction by noticing that $\phi(x_{k_j})\to \phi(\ell)$ and $\phi(x_{k_j})\ge p_{k_j}-1 \to +\infty$.
			
			Thus for all integers $n\in \N^{*}$ there exists $k_n\in \N^{*}$ such that $\forall k\ge k_n, \frac{a_k}{m}\ge n$. We can choose $k_1=1$ and arrange so that the sequence $(k_n)_{n\in \N^{*}}$ is increasing. Let $\chi$ be a smooth increasing convex function such that $\chi(0)=0$ and $\chi(nm)\ge p_{k_{n+1}} \ \forall n\in \N^{*}$. We claim that $\chi\circ \psi > \phi$ on $V'\setminus U$. Indeed, for all $k\in \N^{*}$, we can find some $n\in \N^{*}$ such that $k_n\le k < k_{n+1}$. Then because $k\ge k_n$, we find $a_k\ge nm$ so $\forall x\in A_{p_k}, \psi(x)\ge a_k \ge nm$. We deduce that\begin{equation*}
				\forall x\in A_{p_k}, \chi\circ\psi(x)\ge \chi(nm) \ge p_{k_{n+1}} > \phi(x). \qedhere
			\end{equation*}
		
		\end{proof}

		\begin{proof}[Proof of Proposition \ref{lem-PSH3}-2]
			By assumption, there is an open neighborhood $V$ of $Z$ and a smooth strictly plurisubharmonic function $\phi$ on $\widetilde{V}:=\pi^{-1}(V)$ which is an exhaustion relative to $\widetilde{X}$. Let $U$ be an open neighbourhood of $Z$ whose closure is included in $V$, and $\psi: \widetilde{X} \to \R$ be the function given by the first point of the Proposition.	Write $\widetilde{Z}:=\pi^{-1}(Z)$. 
			
			After replacing $\widetilde{V}$ by a smaller open subset $\widetilde{V'}$, let $\chi$ be an increasing smooth convex function, given by Lemma \ref{lem-PSH4}, such that $\chi(0)=0$ and $\chi\circ\psi > \phi+2$ on $\widetilde{V'}\setminus\widetilde{U}$. Set $\phi_2:=\phi+1$ and \linebreak $\psi_2:=\chi\circ\psi$. Setting $\eta=\frac{1}{2}$, let $M_{(\eta,\eta)}$ be the regularized maximum function defined in \cite[{Lemma I.5.18}]{demaillyComplexAnalyticDifferential}. Recall that this is a smooth symmetric function on $\R^{2}$ with the following properties: $\max(x,y)\le M_{(\eta,\eta)}(x,y)$, and $M_{(\eta,\eta)}(x,y)=y$ whenever $y\ge x+2\eta$. Also, $M_{(\eta,\eta)}(u_1,u_2)$ is (strictly) plurisubharmonic whenever $u_1,u_2$ are (strictly) plurisubharmonic functions. Now let $\lambda:\widetilde{X}\to \R$ be the function defined by\begin{align*}
				\lambda(x) := \left\{ \begin{array}{ll}
					M_{(\eta,\eta)}(\phi_2(x),\psi_2(x)) & \text{ if } x\in \widetilde{V'}, \\
					\psi_2(x) &\text{ if } x \in \widetilde{X}\setminus \widetilde{V}'.
				\end{array}\right.
			\end{align*}
			Notice that $\psi_2\ge\phi_2+2\eta$ on $\widetilde{V'}\setminus\widetilde{U}$ and $\phi_2\ge\psi_2+2\eta$ on some neighbourhood of $\widetilde{Z}$. Hence $\lambda$ is a smooth strictly plurisubharmonic function. For any real number $C$, we have\begin{equation*}
				\{x\in \widetilde{X} \ \vert \ \lambda(x)\le C\} \subset \{x\in \widetilde{X}\setminus\widetilde{U} \ \vert \ \psi_2(x)\le C\} \cup \{x\in \widetilde{V} \ \vert \ \phi_2(x)\le C\}.
			\end{equation*}
			The right-hand side being compact, we infer that $\lambda$ is an exhaustion, which implies that $\widetilde{X}$ is a Stein manifold.
		\end{proof}
		
		\subsection{A strictly plurisubharmonic exhaustion function on the universal cover of $X_\Gamma$}\label{ss-section-univ_Stein}
		From now on, we fix a non-uniform torsion-free lattice $\Gamma$ of $\PU(n,1)$ satisfying the following properties: the space $\HNC/\Gamma$ admits a toroidal compactification $X_\Gamma$ with $\pi_1$-injective boundary tori, the Albanese map $A_{\Gamma}: X_\Gamma \to \Alb(X_\Gamma)$ is an immersion on the open set $\HNC/\Gamma$ and for any boundary torus $T$ of $X_{\Gamma}$, $A_{\Gamma\,{\vert}T}$ is either an immersion or a constant. We will show that for such a lattice $\Gamma$, the universal cover $\widetilde{X_\Gamma}$ of $X_{\Gamma}$ is Stein.\newline
		
		This statement implies Theorem \ref{thm_4_rev_Stein}. Indeed, in Hummel and Schroeder's article, it is noticed that boundary tori are totally geodesic in $X_\Gamma$ for the metric of non-positive Riemannian curvature constructed there, and they are thus $\pi_1$-injected. In combination with Theorem \ref{thm_3_imm_alb}, we deduce that for any torsion-free arithmetic lattice $\Gamma_0$ of $\PU(n,1)$, there exists a finite index subgroup $\Gamma' <\Gamma_0$ all of whose finite index subgroups satisfy these properties. We will also use that toroidal compactifications are projective manifolds.\newline
		
		Here is the global strategy of our proof. Let us denote by $\pi:\widetilde{X_{\Gamma}}\to X_{\Gamma}$ the universal cover of $X_{\Gamma}$. We can construct holomorphic functions on $\widetilde{X_{\Gamma}}$ of the form $f := gh$ where $h$ is a meromorphic function obtained by the theorem of holomorphic convexity with respect to a line bundle \cite[Corollary 4.3]{napierConvexityPropertiesCoverings1990} and $g$ is a holomorphic function which kills the poles of $h$, obtained by post-composing the lift of the Albanese map of $X_{\Gamma}$ by a holomorphic function from $\C^{b_1(X_{\Gamma})/2}$ to $\C$. This strategy gives enough holomorphic functions to construct a strictly plurisubharmonic exhaustion $\psi$ on $\widetilde{X_{\Gamma}}\setminus\pi^{-1}(U)$ where $U$ is an open neighborhood of the set $Z$ consisting of all points $x \in X_{\Gamma}$ whose fiber $A_{\Gamma}^{-1}(A_{\Gamma}(x))$ is not finite (Proposition \ref{lem-PSH3}-1). We also construct by hand a strictly plurisubharmonic exhaustion function $\phi$ on the lift of an open neighborhood of $Z$ (Lemmas \ref{lem-PSH1}, \ref{lem-PSH2}) and finally we glue $\phi$ and $\psi$ (Lemma \ref{lem-PSH4}) to obtain a strictly plurisubharmonic exhaustion function on $\widetilde{X_\Gamma}$ (Proposition \ref{lem-PSH3}-2). We use the fact that $A_{\Gamma}$ restricted to any boundary torus $T$ is either an immersion or a constant in order to get a nice description of $Z$ (Lemma \ref{lem-PSH2} and its proof).\newline
		
		To prove Theorem \ref{thm_4_rev_Stein}, we show that the Albanese map of $X_{\Gamma}$ satisfies all the hypotheses of Proposition \ref{lem-PSH3}. This is the content of Lemma
		\ref{lem-PSH2}, whose proof uses Lemmas \ref{lem-PSH1}. Lemma \ref{lem-PSH1} does not use the assumptions on the Albanese map which are described in the beginning of the Subsection, contrary to Lemma \ref{lem-PSH2}. 
				
		\begin{lemma}\label{lem-PSH1}
			Let $T$ be a boundary torus of $X_\Gamma$. If $V$ is a sufficiently small open neighbourhood of $T$, its preimage $\widetilde{V}:=\pi^{-1}(V)$ is a disjoint union of neighbourhoods of (copies of) $\C^{n-1}$. We can choose $V$ so that each of these open neighbourhood admits a strictly plurisubharmonic function which is of exhaustion relative to $\widetilde{X_\Gamma}$.
		\end{lemma}
		
		 For Lemma \ref{lem-PSH2}, we recall that the definition of a relative exhaustion has been given in Subsection \ref{sssection_res_gen}.
		
		\begin{lemma}\label{lem-PSH2}
		 Define\begin{equation*}
				Z := \{x\in X_{\Gamma} \ \vert \ A_{\Gamma}^{-1}(A_{\Gamma}(x)) \text{ is not finite}\}.
			\end{equation*}
			Then the set $Z$ is a finite number of fibers of the map $A_{\Gamma}$, and it consists of some boundary tori together with a finite number of points. In particular, there exists a neighbourhood $V$ of $Z$ as well as a smooth strictly plurisubharmonic function $\phi$ on $\widetilde{V}:=\pi^{-1}(V)$ which is an exhaustion relative to $\widetilde{X_{\Gamma}}$.
		\end{lemma}

		\begin{proof}[Proof of Lemma \ref{lem-PSH1}]
			By hypothesis, for any boundary torus $T$ of $X_{\Gamma}$, the map $\pi_1(T)\to \pi_1(X_{\Gamma})$ induced by the inclusion $T\hookrightarrow X_{\Gamma}$ is injective, therefore the preimage of $T$ in the universal cover $\widetilde{X_{\Gamma}}$ of $X_{\Gamma}$ is a disjoint union of copies of $\C^{n-1}$.

			Let $p:\C^{n-1}\to T$ be the universal cover of a boundary torus $T$ of $X_\Gamma$, $V$ a neighbourhood of $T$ identified with the unit disk bundle of a holomorphic line bundle on $T$ endowed with a Hermitian metric $h$ as in Proposition \ref{prop_toroidal_compactification}, and $F:V \to T$ the map defining this bundle. Then we have the following commutative diagram:
			\begin{equation*}
				\begin{tikzcd}
					\widetilde{V}\arrow{d}{\widetilde{F}}\arrow{r}{\pi}&V\arrow{d}{F} \\
					\C^{n-1} \arrow{r}{p}  & T,
				\end{tikzcd}
			\end{equation*}
			where $\pi:\widetilde{V}\to V$ is the universal cover of $V$ and $\widetilde{F}$ the lift of $F$. We identify $\widetilde{V}$ with a connected component of the preimage of $V$ in $\widetilde{X_\Gamma}$ and $\C^{n-1}$ with the connected component of the preimage of $T$ in $\widetilde{X_\Gamma}$ contained in $\widetilde{V}$. We will show that some neighborhood of the zero section in $\widetilde{V}$ has a strictly plurisubharmonic function $\phi$ which is an exhaution relatively to $\widetilde{X_\Gamma}$. This function $\phi$ will be the sum of the pullback by $\widetilde{F}$ of a strictly plurisubharmonic exhaustion function on $\C^{n-1}$, for instance $z\mapsto \lVert z\rVert^{2}$, and a bounded function on the fibers, defined to be the pullback by $\pi$ of the squared hermitian norm on $V$, which will be strictly plurisubharmonic in a neighborhood of the zero section. Explicitely, for $\widetilde{v}\in \widetilde{V}$, let us put $\phi(\widetilde{v}) := \lVert \widetilde{F}(\widetilde{v}) \rVert^{2} + N(\pi(\widetilde{v}))$, where $N:V\to \R_+$ is the function $v\mapsto h(v,v)$. Then the function $\phi$ is an exhaustion relative to $\widetilde{X_{\Gamma}}$.
			
			We now assert that $\phi$ is strictly plurisubharmonic on $\C^{n-1}$. Indeed, working in coordinates, this amount to saying that the map\begin{equation*}
				(a,v) \in \C\times\C^{n-1} \longmapsto \lVert v\rVert^{2} + \exp({h(v)})\lvert a\rvert^{2}
			\end{equation*}
			is strictly plurisubharmonic at all points of the form $(0,v)$, which is easily checked. Since $i\partial\bar\partial \phi$ is invariant under the action of $\pi_1(T)$ on $\widetilde{V}$, which is co-compact on some closed neighbourhood of $\C^{n-1}$, we deduce that $\phi$ is strictly plurisubharmonic on a neighbourhood of $\C^{n-1}$ of the form $\pi^{-1}(V')$ for some open subset $V'\subset V$.
		\end{proof}
		
		\begin{proof}[Proof of Lemma \ref{lem-PSH2}]
			Let $T_1,\dots,T_k$ be the boundary tori of $X_{\Gamma}$ on which the Albanese map $A_{\Gamma}$ is constant, of values $a_1,\dots,a_k$. We claim that $Z=A_{\Gamma}^{-1}(\{a_1,\dots,a_k\})$, which implies that it is an analytic set. Indeed, if $x\in Z$, there exists a smooth immersed curve $z\mapsto \gamma(z)$ in $A_{\Gamma}^{-1}(A_{\Gamma}(x))$. For all $z$, $dA_{\Gamma}\cdot\gamma'(z)=0$ so $\gamma$ is included in a torus. The Albanese map is an immersion in restriction to all tori except $T_1,\dots,T_k$ so $\gamma$ must be included in one of the $T_i$, and $A_\Gamma(x)=a_i$. This proves one inclusion, the other being immediate. In particular, any smooth curve included in $Z$ is either constant or included in one of the tori $T_1,\dots,T_k$ on which the Albanese map $A_{\Gamma}$ is constant. We deduce that $Z$ is the union of the tori $T_1,\dots,T_k$ with a discrete subset of $X_{\Gamma}\setminus\bigcup_{i=1}^{k}T_i$. This discrete subset has to be finite, since a compact analytic set has finitely many irreducible components.
			
			The last part of the lemma then follows from Lemma \ref{lem-PSH1} and the fact that each point has a Stein neighborhood.
		\end{proof}

\bigskip
\bigskip
\begin{small}\begin{tabular}{l}
		William, Sarem\\
		Univ. Grenoble Alpes, CNRS, IF, 38000 Grenoble, France\\
		william.sarem at univ-grenoble-alpes.fr
	\end{tabular}
\end{small}		
\end{document}